\documentclass[11pt]{amsart}

\usepackage[margin=1.0in]{geometry}
\usepackage{microtype}

\usepackage{graphicx}
\usepackage{amssymb}
\usepackage{enumerate}
\usepackage{times}
\usepackage{cite}

\usepackage{xcolor}

\usepackage[colorlinks=true, pdfstartview=FitV, linkcolor=blue, citecolor=blue, urlcolor=blue]{hyperref}

\def\Xint#1{\mathchoice
{\XXint\displaystyle\textstyle{#1}}%
{\XXint\textstyle\scriptstyle{#1}}%
{\XXint\scriptstyle\scriptscriptstyle{#1}}%
{\XXint\scriptscriptstyle\scriptscriptstyle{#1}}%
\!\int}
\def\XXint#1#2#3{{\setbox0=\hbox{$#1{#2#3}{\int}$}
\vcenter{\hbox{$#2#3$}}\kern-.5\wd0}}

\def\avgint{\Xint-}

\makeatletter
\newcommand*\rel@kern[1]{\kern#1\dimexpr\macc@kerna}
\newcommand*\widebar[1]{%
  \begingroup
  \def\mathaccent##1##2{%
    \rel@kern{0.8}%
    \overline{\rel@kern{-0.8}\macc@nucleus\rel@kern{0.2}}%
    \rel@kern{-0.2}%
  }%
  \macc@depth\@ne
  \let\math@bgroup\@empty \let\math@egroup\macc@set@skewchar
  \mathsurround\z@ \frozen@everymath{\mathgroup\macc@group\relax}%
  \macc@set@skewchar\relax
  \let\mathaccentV\macc@nested@a
  \macc@nested@a\relax111{#1}%
  \endgroup
}
\makeatother

\newtheorem{theorem}{Theorem}[section]
\newtheorem{lemma}[theorem]{Lemma}

\newtheorem{corollary}[theorem]{Corollary}

\newtheorem{definition}[theorem]{Definition}

\theoremstyle{definition}

\theoremstyle{remark}
\newtheorem{remark}[theorem]{Remark}

\numberwithin{equation}{section}
 \allowdisplaybreaks

\newcommand{\R}{\mathbb R}

\newcommand{\BB}{\mathbf B}

\newcommand{\Q}{\mathcal Q}

\newcommand{\F}{\mathcal F}

\newcommand{\A}{\mathcal A}
\newcommand{\B}{\mathcal B}

\DeclareMathOperator*{\esssup}{ess\,sup}

\DeclareMathOperator{\op}{{op}}

\DeclareMathOperator{\clconv}{\overline{conv}}

\newcommand{\K}{\mathcal{K}}

\newcommand{\Rdf}{\mathcal{R}}

\newcommand{\loc}{\text{loc}}

\title[Sharp constant off-diagonal extrapolation for Muckenhoupt bases]
{Sharp constant off-diagonal extrapolation \\  for Muckenhoupt bases}

\author{David Cruz-Uribe OFS}
\address{Department of Mathematics \\ The University of Alabama \\
Tuscaloosa, AL 35487, USA}
\email{dcruzuribe@ua.edu}

\author{Fatih \c{S}irin}
\address{Department of Mathematics \\ Hali\c{c} University \\ 5. Levent Campus, Istanbul, 34060, T\"urkiye}
\email{fatihsirin@halic.edu.tr}

\thanks{The first author is partially supported by a Simons Foundation Travel Support for Mathematicians Grant and by NSF grant DMS-2349550.  The second author was supported by 2219 - International Postdoctoral Research Fellowship Program for Turkish Citizens from T\"ubitak, the Scientific and Technological Research Council of T\"urkiye.}

\subjclass[2010]{42B25, 42B35}

\keywords{convex-set valued functions,
  maximal operators, off-diagonal inequalities,
  Muckenhoupt weights, matrix weights,
  Rubio de Francia extrapolation, sharp constants}

\begin{document}

\maketitle

\section{Introduction}
In this paper we prove a sharp constant version of off-diagonal extrapolation for matrix weights, improving the result we proved in~\cite{MR4972891}.  To put our result into context, we briefly sketch some recent results.    By a matrix weight $W$ we mean a $d\times d$ self-adjoint matrix function that is finite and invertible almost everywhere:  that is, its coefficients are measurable and finite almost everywhere, and  $W(x)$ is invertible for almost every $x$.    Given $1\leq p \leq\infty$,  define the space $L^p(\R^n, W)$ to be the Banach space of measurable, vector-valued functions $f :\R^n \rightarrow \R^d $  with the norm
\[ \|f\|_{L^p(\R^n,W)} = \|Wf\|_{L^p(\R^n,\R^d)} = \bigg(\int_{\R^n} |W(x)f(x)|^p\,dx \bigg)^{\frac{1}{p}} < \infty.  \]
Note that if $p=\infty$, the first equality still makes sense.  For $1<p<\infty$,  a matrix weight $W$ is in $\A_p$ if
\[ [W]_{\A_p} = \sup_Q \bigg(\avgint_Q \bigg( \avgint_Q
    |W(x)W^{-1}(y)|_{\op}^{p'}\,dy\bigg)^{\frac{p}{p'}}\,dx\bigg)^{\frac{1}{p}} <
    \infty, \]
  where the supremum is taken over all cubes with sides parallel to the coordinate axes.
  Define $\A_1$ and $\A_\infty$  to be matrix weights that satisfy
  \[ [W]_{\A_1} = \sup_Q \esssup_{x\in Q} \avgint_Q
    |W^{-1}(x)W(y)|_{\op}\,dy < \infty,  \]
and
\[ [W]_{\A_\infty} =\sup_Q \esssup_{x\in Q} \avgint_Q
  |W(x)W^{-1}(y)|_{\op}\,dy < \infty.  \]

In~\cite{bownik-cruz-uribe}, the first author and Bownik resolved an open problem due to Nazarov and Treil~\cite{MR1428988}, proving the following generalization of Rubio de Francia extrapolation.

\begin{theorem} \label{thm:matrix-extrapolation-intro} 
Let $T$ be a  linear operator.  Suppose that for some $p_0$,
  $1 \leq p_0\leq\infty$, there exists an increasing function
  $N_{p_0}$ such that for every $W_0\in \A_{p_0}$,
\begin{equation}\label{ext1}
 \|Tf\|_{L^{p_0}(\R^n, W_0)}
  \leq N_{p_0}([W_0]_{\A_{p_0}})\|f\|_{L^{p_0}(\R^n, W_0)}. 
  \end{equation}
  Then for all $p$, $1<p<\infty$, and 
  for all $W\in \A_p$,
\begin{equation}\label{ext2}
    \|Tf\|_{L^p(\R^n,W)}
  \leq N_{p}(n,d,p_0,[W]_{\A_{p}})\|f\|_{L^{p}(\R^n,W)},
  \end{equation}
  where
  \[ N_p (n,d,p_0,[W]_{\A_{p}})
    = C(p,p_0)N_{p_0}\bigg(C(n,d,p,p_0)[W]_{\A_p}^{\max\big\{\frac{p}{p_0},\frac{p'}{p_0'}\big\}}\bigg). \]
\end{theorem}

A key feature of this result is that the constant $N_p$ depends on the power of the $\A_p$ characteristic $[W]_{\A_{p}}$ with the same exponent that appears in the scalar case.  In the scalar case this is referred to as ``sharp constant" extrapolation, since if one starts with the best possible dependence on $[w]_{A_2}$ in $L^2(w)$ estimates for singular integrals, then this theorem yields the best possible dependence on $[w]_{A_p}$ in the $L^p(\R^n,w)$ estimates, $1<p<\infty$.  This fact played a very important role in the first proofs of the so-called $A_2$ conjecture:  see Hyt\"onen~\cite{hytonenP2010,MR3204859}.  

In~\cite{MR4972891}, we extended Theorem~\ref{thm:matrix-extrapolation-intro} to off-diagonal extrapolation, generalizing the corresponding scalar theorem first proved by Harboure, Mac\'\i as, and Segovia~\cite{harboure-macias-segovia88}.  

\begin{theorem} \label{thm:off-diag-extrapol-cubes}
Let $T$ be a linear operator.   Fix $s$, $1\leq s<\infty$.  Suppose that for some $ p_0, q_0, 1 \leq p_0 \leq q_0 \leq \infty $, where $\frac{1}{p_0}-
\frac{1}{q_0}=\frac{1}{s'}$ and $r_0=\frac{q_0}{s}$, there exists a positive increasing function $ N_{p_0,q_0} $ such that for every matrix weight $W_0$ with $ W_0^s \in \A_{r_0} $,
\begin{equation} \label{eqn:odec1}
\|Tf\|_{L^{q_0}(\R^,W_0)} 
 \leq N_{p_0, q_0} \left( [W_0^s]_{\A_{r_0}}\right)
\|f\|_{L^{p_0}(\R^,W_0)} 
\end{equation}
Then for all $ p,\, q $ such that $ 1 < p \leq q < \infty $, $ \frac{1}{p} - \frac{1}{q} = \frac{1}{p_0} - \frac{1}{q_0}$, $r=\frac{q}{s}$, and for all matrix weights $W$ with $ W^s \in \A_{r} $,
\begin{equation} \label{eqn:odec2}
\|Tf\|_{L^{q}(\R^,W)} 
\leq N_{p, q}(n,d,p_0,q_0,[W^s]_{\A_{r}})
\|f\|_{L^{p}(\R^,W)} 
\end{equation}
where 
\[
N_{p, q}(n,d,p_0,q_0,[W^s]_{\A_{r}}) = C(d,p,q,p_0,q_0)
N_{p_0, q_0}\big(C(n,d,p,q,p_0, q_0) [W^s]_{\A_{r}}^{\frac{r}{r_0} + \frac{r'}{r'_0}}\big).
\]
\end{theorem}

We call attention to two features of this proof.  First, in the scalar case, the condition $[w^s]_{\A_r}$ is equivalent to the $\A_{p,q}$ condition introduced by Muckenhoupt and Wheeden.  But in the matrix case, it was shown in~\cite[Example~2.8]{MR4972891} that the corresponding condition matrix $\A_{p,q}$ condition,
\[
[W]_{\A_{p,q}} = \sup_{Q} \left( \avgint_Q \left( \avgint_{Q} |W(x) W^{-1}(y)|^{p'}_{\op} \,dy \right)^{\frac{q}{p'}} \,dx \right)^{\frac{1}{q}}  < \infty,
\]
is strictly weaker than the condition that $W^s\in \A_r$.  It is an open problem whether Theorem~\ref{thm:off-diag-extrapol-cubes} is true with this weaker assumption.  

Second, and more importantly for the purposes of this paper, in Theorem~\ref{thm:off-diag-extrapol-cubes} the constant $N_{p,q}$ is larger than the sharp constant found in the scalar case by  Lacey, Moen, P\'erez, and Torres~\cite{MR2652182}.  They proved that the exponent on $[w^s]_{\A_r}$ in the sharp constant is $\max \big\{ \frac{r}{r_0},\frac{r'}{r'_0} \big\}$.    Our  main result is showing that this exponent is also gotten in the matrix case, improving Theorem~\ref{thm:off-diag-extrapol-cubes} to a sharp constant result.  

\begin{theorem} \label{thm:off-diag-extrapol-sharp-cubes}
Let $T$ be a sublinear operator.   Fix $s$, $1\leq s<\infty$.  Suppose that for some $ p_0, q_0, 1 \leq p_0 \leq q_0 \leq \infty $, where $\frac{1}{p_0}-
\frac{1}{q_0}=\frac{1}{s'}$ and $r_0=\frac{q_0}{s}$, there exists a positive increasing function $ N_{p_0,q_0} $ such that for every matrix weight $W_0$ with $ W_0^s \in \A_{r_0} $, inequality~\eqref{eqn:odec1} holds.  Then for all $ p,\, q $ such that $ 1 < p \leq q < \infty $, $ \frac{1}{p} - \frac{1}{q} = \frac{1}{p_0} - \frac{1}{q_0}$, $r=\frac{q}{s}$, and for all matrix weights $W$ with $ W^s \in \A_{r} $, inequality~\eqref{eqn:odec2} holds, with 
\[ N_{p,q}(n,d,p_0,q_0, [W^s]_{\A_{r,\B}}) 
    = C(p,q,p_0,q_0) N_{p_0,q_0}\big(C(n,d,p,q,p_0,q_0) 
    [W^s]_{\A_{r}}^{\max \big\{ \frac{r}{r_0},\frac{r'}{r'_0} \big\}} \big). \]
\end{theorem}

\medskip

We actually prove a more general result.  First, in~\cite{MR4972891}, following~\cite{MR2797562} in the scalar case,  we introduced the concept of a matrix Muckenhoupt basis and proved our results in this general context.  We do the same here.  By a basis $\B$ we mean a collection of open sets.  Given a basis $\B$ and a matrix weight $W$,  define the Christ-Goldberg maximal operator 
\[ M_{W,\B}f(x) = \sup_{B\in \B} \avgint_B |W(x)W^{-1}(y)f(y)|\,dy \cdot \chi_B(x). \]
Define the  class of $\A_{p,\B}$ weights
by replacing cubes with sets $B\in \B$ in the above definitions: e.g., for $1<p<\infty$,
\[ [W]_{\A_{p,\B}} = \sup_{B\in \B} \bigg(\avgint_B \bigg( \avgint_B
    |W(x)W^{-1}(y)|_{\op}^{p'}\,dy\bigg)^{\frac{p}{p'}}\,dx\bigg)^{\frac{1}{p}} <
    \infty. \]

\begin{definition} \label{defn:matrix-muckenhoupt}
   A basis $\B$ is a matrix Muckenhoupt basis if for every $1<p<\infty$ and matrix weight $W\in \A_{p,\B}$, $M_{W,\B} : L^p(\R^n,\R^d) \rightarrow L^p(\R^n)$, and there exists an increasing function $K_{p,\B}$ such that 
    \[\bigg(\int_{\R^n}|M_{W,\B}f(x)|^p\,dx \bigg)^{\frac{1}{p}}
    \leq 
    K_{p,\B}([W]_{\A_{p,\B}})\bigg(\int_{\R^n}|f(x)|^p\,dx\bigg)^{\frac{1}{p}}.
    \]
\end{definition}

\begin{remark}
    As we showed in~\cite{MR4972891}, beside the basis of cubes, the basis of dyadic cubes and the basis of rectangles (i.e., parallelpipeds) in $\R^n$ with sides parallel to the coordinate axes, are bases.  It is an interesting open question to prove that the Zygmund basis in $\R^3$, consisting of parallelpipeds with sides of length $s,\,t,\,st$, $s,\,t>0$, is a matrix Muckenhoupt basis.
\end{remark}

Second, we prove our results using the abstract concept of a family $\F$ of extrapolation pairs.  These were used in~\cite{MR2797562} for scalar-weighted extrapolation, and in~\cite{bownik-cruz-uribe, MR4972891} for matrix-weighted extrapolation.  Hereafter, $\F$ will denote a family of pairs $(f,g)$ of
measurable, vector-valued functions.  If we write an
inequality of the form
\[ \|f\|_{L^q(\R^n,W)} \leq C\|g\|_{L^p(\R^n,W)}, \qquad (f,g) \in
    \F, \]
  we mean that this inequality holds for all pairs $(f,g)\in\F$ such
  that 
  \[ 0 < \|f\|_{L^q(\R^n,W)}, \|g\|_{L^p(\R^n,W)} < \infty. \]
  The constant,
  whether given explicitly or implicitly, is assumed to be independent
  of the pair $(f,g)$ and to depend only on $[W]_{\A_p}$ and not on
  the particular matrix weight $W$. 

\begin{theorem} \label{thm:sharp-off-diag-extrapol-pairs}
Let $\B$ be a matrix Muckenhoupt basis and let $\F$ be a family of extrapolation pairs.   Fix $s$, $1\leq s<\infty$. Suppose that for some $ p_0, q_0, 1 \leq p_0 \leq q_0 \leq \infty $, where $\frac{1}{p_0}-\frac{1}{q_0}= \frac{1}{s'}$ and $r_0=\frac{q_0}{s}$, there exists a positive increasing function $ N_{p_0,q_0}$ such that for every matrix weight $W$ with $ W_0^s \in \A_{r_0,\B} $,
\begin{equation} \label{eqn:sharp-extrapol0}
\|f\|_{L^{q_0}(\R^,W_0)} 
\leq N_{p_0, q_0} \left( [W_0^s]_{\A_{r_0,\B}}\right)
\|g\|_{L^{p_0}(\R^,W_0)}, \quad (f,g) \in \F.
\end{equation}
Then for all $ p,\, q $ such that $ 1 < p \leq q < \infty $, $ \frac{1}{p} - \frac{1}{q} = \frac{1}{p_0} - \frac{1}{q_0}$, $r=\frac{q}{s}$, and for all matrix weights $W$ with $ W^s \in \A_{r,\B} $,
\begin{equation} \label{eqn:Sharp-extrapol1}
\|f\|_{L^{q}(\R^,W)} 
\leq N_{p, q}(\B,p_0,q_0,[W^s]_{\A_{r,\B}})
\|g\|_{L^{p}(\R^n,W)}, \qquad (f,g) \in \F,
\end{equation}
where if $p<p_0$,
\begin{equation*} 
    N_{p,q}(\B,p_0,q_0, [W^s]_{\A_{r,\B}}) 
    = 2^{\left( \frac{1}{s}+1 \right)\frac{q}{p} } N_{p_0,q_0}\big(C(d,r) [W^s]_{\A_{r,\B}}^{\frac{r}{r_0}}  
    K_{r,\B}( [W^s]_{\A_{r,\B}})^{1-\frac{r}{r_0}}\big), 
\end{equation*}
and if $p>p_0$,
\begin{equation*}
    N_{p,q}(\B,p_0,q_0, [W^s]_{\A_{r,\B}}) 
    = 2^{r'\frac{p-p_0}{pp_0}} N_{p_0,q_0}\big( C(d,r) [W^s]_{\A_{r,\B}}^{\frac{r'}{r'_0}}  
    K_{r',\B}([W^{-s}]_{\A_{r',\B}}) )^{1-\frac{r'}{r'_0}} \big)
\end{equation*}
\end{theorem}

\begin{remark}
    In the scalar case, it is straightforward to use extrapolation pairs to prove results for operators via an approximation argument:  see~\cite{preprint-DCU, MR2797562} for details.  For matrix weights there are some subtleties, especially for nonlinear operators.  We will discuss these briefly below when we prove Theorem~\ref{thm:off-diag-extrapol-cubes} as a consequence of Theorem~\ref{thm:sharp-off-diag-extrapol-pairs}.  For a detailed discussion, see~\cite[Section~10]{bownik-cruz-uribe}.  For an alternative approach that avoids extrapolation pairs, see Nieraeth~\cite{MR4966577}.
\end{remark}

\medskip

The remainder of this paper is organized as follows.  In Section~\ref{section:matrix} we prove Theorem~\ref{thm:sharp-off-diag-extrapol-pairs} and then derive Theorem~\ref{thm:off-diag-extrapol-cubes} as a special case.  The basic idea of the proof is to adapt the diagonal (that is, when $p=q$) sharp constant proof in~\cite{bownik-cruz-uribe} to the off-diagonal case.  That proof, in turn, was gotten by adapting the proof of sharp constant extrapolation in the scalar case given in~\cite[Section~3.4]{MR2797562}.   The original proof of sharp constant, off-diagonal extrapolation in~\cite{MR2652182} relies on a very technical lemma about Muckenhoupt $A_p$ weights, which was not easily adapted to the matrix case.  Therefore, in order to prove our desired result, we first had to extend the proof of sharp constant extrapolation in~\cite[Section~3.4]{MR2797562} to the off-diagonal case.  We believe that this result is of independent interest, so we have included it in Appendix~\ref{section:scalar}.  Readers who want an overview of the main ideas of the the proof of Theorem~\ref{thm:sharp-off-diag-extrapol-pairs} and an understanding of where  the technical obstacles  appear in the matrix-weight setting, should first read the proof of Theorem~\ref{thm:sharp_constant-scalar} and Remark~\ref{remark:scalar-matrix} following the proof.

\begin{remark}
    Duoandikoetxea~\cite{10.1016/j.jfa.2010.12.015} gave another proof of sharp constant, off-diagonal extrapolation, and extended it to the case $q<p$.  It is an interesting problem to see if his ideas can be extended to the matrix setting.
\end{remark}

Throughout this paper, we use the notation from~\cite{bownik-cruz-uribe,MR4972891}; we will also refer the reader to these papers for a number of technical results needed in our proofs.
 In Euclidean space, the constant $n$ will denote
the dimension of $\R^n$, which will be the domain of our functions.  The value
$d$ will denote the dimension of vector and set-valued functions.     For $1\leq p \leq \infty$, $L^p(\R^n)$ will denote the Lebesgue space
 of scalar functions, and $L^p(\R^n,\R^d)$ will denote the Lebesgue
 space of vector-valued functions.  By a scalar weight we mean a measurable function $w$ such that $0<w(x)<\infty$ almost everywhere.  We define the space $L^p(\R^n,w)$ to consist of all measurable functions $f$ such that $wf\in L^p(\R^n)$.  Note that this is not the standard definition, which treats weights as measures rather than multipliers,
but it brings our scalar definition in line with our definition of the matrix weighted norm.

 Given $v=(v_1,\ldots,v_d)^t \in \R^d$, the Euclidean norm of $v$ will
 be denoted by $|v|$.  The standard orthonormal basis in $\R^d$ will
 be denoted by $\{e_i\}_{i=1}^d$.  The open unit ball 
 $\{ v \in \R^d : |v|<1\}$ will be denoted by ${\mathbf B}$ and its
 closure by $\overline{\mathbf{B}}$.  Matrices will be $d\times d$
 matrices with real-valued entries.  
 Given
 two quantities $A$ and $B$, we will write $A \lesssim B$, or
 $B\gtrsim A$ if there is a constant $c>0$ such that $A\leq cB$.  If
 $A\lesssim B$ and $B\lesssim A$, we will write $A\approx B$.

\section{Off-diagonal matrix extrapolation:  proof of Theorem~\ref{thm:sharp-off-diag-extrapol-pairs}}
\label{section:matrix}

In this section we  first prove Theorem~\ref{thm:sharp-off-diag-extrapol-pairs}; afterwards, we   derive Theorem~\ref{thm:off-diag-extrapol-sharp-cubes} as an immediate consequence.  For the proof we need a handful of definitions and results from~\cite{bownik-cruz-uribe}.  We will give a precise statement of them here, before giving the actual proof.

We begin with some basic definitions.  The collection of all convex, symmetric sets in $\R^d$ is denoted by $\K_{cs}(\R^d)$.  The subset of bounded sets in  $\K_{cs}(\R^d)$  is denoted by $\K_{bcs}(\R^d)$.  By a measurable convex-set valued function we mean a map $F : \R^n \to \K_{cs}(\R^d)$ such that given any open set $U\in \R^d$, $F^{-1}(U)= \{ x\in \R^n : F(x)\cap U \neq \emptyset\}$ is Lebesgue measurable. It is integrably bounded if there exists $k\in L^1(\R^n)$ such that $F(x)\subset k(x)\overline{\BB}$ for a.e.~$x$. If $k$ is in $L^1_{\loc}$, we say $F$ is locally integrably bounded. Given a matrix weight $W$ and $1\leq p\leq \infty$, $L^p_\K(\R^n,W)$ consists of all measurable convex-set valued functions $F$ such that 
\[ \|F\|_{L^p_\K(\R^n, W) } = \||WF|\|_{L^p(\R^n)} = \bigg(\int_{\R^n} |W(x)F(x)|^p\,dx\bigg)^{\frac{1}{p}} < \infty, 
\]
where $|W(x)F(x)| = \sup\{ |W(x)v| : v\in F(x)\}$.   If $W=I_d$,  we write $L^p_\K(\R^n)=L^p_\K(\R^n,I_d)$.  

We say that a convex-set valued function $F$ is ball-valued if there exists a non-negative, measurable scalar function $f$ such that $F=f\overline{\BB}$.  More generally, $F$ is ellipsoid-valued if there exists a a matrix weight $W$ such that $F=fW\overline{\BB}$.  

We now consider operators mapping meausurable convex-set valued functions to themselves.  Such an operator $T$ is sublinear if for any two such functions $F$ and $G$, $T(F+g)(x) \subset TF(x) + TG(x)$, and $T(\alpha F)(x)= \alpha TF(x)$ where $\alpha \in \R$ and the sum is the Minkowski sum of sets. It is monotone if, whenever $F(x) \subset G(x)$, $TF(x)\subset TG(x)$.

Our principal such operator is the convex-set valued maximal operator.  Recall that the Aumann integral of a locally integrably bounded function $F : \R^n \to \K_{bcs}$ on a bounded set $\Omega$ is defined to be the convex set
\[  \int_\Omega F(y)\,dy =  \bigg\{ \int_\Omega f(y)\,dy : f \in S^1(\Omega) \bigg\}; \]
here, $S^1(\Omega)$ is the set of vector-valued functions $f\in L^1\Omega)$ such that $f(y)\in F(y)$, $y\in \Omega$. 
Given a basis $\B$, define the convex-set valued maximal operator $M_\B$ to be
\[ M_\B F(x) = \clconv\bigg(\bigcup_{B\in \B} F(y)\,dy \cdot \chi_B(x)\bigg), \]
where $\clconv(A)$, $A\subset \R^d$, is the closed, convex hull of $A$. As proved in~\cite[Section~3]{MR4972891}, $M_\B$ is a measurable convex-set valued function that is sublinear and monotone.  Moreover, we have the following norm inequality.

\begin{theorem}\label{thm_6.9}
  Let $\B$ be a Muckenhoupt basis.  Fix $1<p\leq \infty$   Given a matrix $ W \in A_{p,\B} $, $ 1 < p < \infty $, the convex-set valued maximal operator $ M_\B: L^p_{\K}(\R^n, W) \to L^p_{\K}(\R^n, W) $ is bounded:  
  \[ \|M_\B F\|_{L^p_{\K}(\R^n, W)} \leq C(d)K_{p,\B}([W]_{\A_p})\|MF\|_{L^p_{\K}(\R^n, W)}, \]
  where $K_{p,\B}$ is the constant from Definition~\ref{defn:matrix-muckenhoupt}.
  Moreover, if $\B=\Q$, the basis of cubes, then
  \[ \|M_\Q F\|_{L^p_{\K}(\R^n, W)} \leq C(n,d,p)[W]_{\A_p}^{p'}\|MF\|_{L^p_{\K}(\R^n, W)}. \]
\end{theorem}

As part of the extrapolation argument, we will need to pass from the maximal operator to a closely related function that is ball or ellipsoid-valued.  
Given a matrix weight $W$, we define the ellipsoid-valued operator $N_W$ by $N_WF(x) = |W(x)F(x)|W^{-1}(x)\overline{\BB}$.  Then $N_W$ is sublinear and monotone, and an isometry: for all $p$, $\|N_WF\|_{L^p_\K(\R^n,W)} = \|F\|_{L^p_\K(\R^n,W)}$.  See~\cite[Lemma~8.4]{bownik-cruz-uribe}.  

We will also need the following convex-set valued version of the Rubio de Francia iteration operator.  The following result was proved in~\cite[Theorem~7.6, Lemma~8.5]{bownik-cruz-uribe}.

\begin{theorem}\label{thm_7.6}
   Given a matrix weight and  $1 \leq p < \infty$, suppose $T$ is a convex-set valued operator with the following properties:
\begin{enumerate}
    \item $T : L^p_\K(\R^n, W) \to L^p_\K(\R^n, W)$ with norm $\|T\|_{L^p_\K(\R^n, W)}$;
    \item $T$ is sublinear and monotone.
\end{enumerate}

Given $G \in L^p_\K(\R^n, G)$, define
\begin{equation}
    SG(x) = \sum_{k=0}^\infty 2^{-k} \|T\|_{L^p_\K(\R^n, W)}^{-k} T^k G(x),
\end{equation}
where $T^k G = T \circ T \circ \cdots \circ T G$ for $k \geq 1$ and $T^0 G(x) = G(x)$. Then this series converges in $L^p_\K(\R^n, W)$ and $SG : \R^n \to \K_{bcs}(\R^d)$ is a measurable mapping. Moreover, it has the following properties:
\begin{enumerate}
    \item $G(x) \subset SG(x)$;
    \item $\|SG\|_{L^p_\K(\R^n, \rho)} \leq 2 \|G\|_{L^p_\K(\R^n, W)}$;
    \item $T(SG)(x) \subset 2 \|T\|_{L^p_\K(\R^n, W)} SG(x)$.
\end{enumerate}
Further, if $TG$ is an ellipsoid or ball-valued function whenever $G$ is, then $SG$ is also ellipsoid or ball-valued.
\end{theorem}

Finally, we need three results about matrix $\A_p$ weights.  All of these results were proved for the basis cubes $\Q$, but the same proof holds for any Muckenhoupt basis $\B$.
The first was proved in \cite[Corollary~8.9]{bownik-cruz-uribe}.

\begin{lemma} \label{matrix-weightsA_r}
    Given a Muckenhoupt basis $\B$, $1<p<\infty$ and $W\in \A_{p,\B}$, suppose there is a scaler function $s$ such that $W_1=sW\in \A_{\infty,\B}$. Then for $p<p_0<\infty$, $\bar{W}=W^\frac{p}{p_0}W_1^{1-\frac{p}{p_0}}\in \A_{p_0,\B}$; moreover,
    \[
    [\bar{W}]_{\A_{p_0,\B}} \le C(d) [W]_{\A_{p,\B}}^{\frac{p}{p_0}} [W_1]_{\A_{\infty,\B}}^{1-\frac{p}{p_0}}.
    \]
    Similarly, if there exists a scalar function $r$ such that $W_0=rW\in \A_{1,\B}$, then for every $1<p_0<p$, $\bar{W}=W_0^{1-\frac{p'}{p'_0}} W^\frac{p'}{p'_0} \in \A_{p_0,\B}$; moreover,
    \[
    [\bar{W}]_{\A_{p_0,\B}} \le C(d)  [W_0]_{\A_{1,\B}}^{1-\frac{p'}{p'_0}} [W]_{\A_{p,\B}}^{\frac{p'}{p'_0}}.
    \]
\end{lemma}

The second is a well-known property of matrix weights; see~\cite[Proposition~6.4]{bownik-cruz-uribe}.

\begin{lemma} \label{cor:duality}
    Given a basis $\B$ and $1\leq p \leq \infty$, a matrix weight $W$ satisfies $W\in \A_{p,\B}$ if and only if $W^{-1}\in \A_{p',\B}$ and in this case $[W]_{\A_{p,\B}}\approx [W^{-1}]_{\A_{p',\B}}$.
\end{lemma}

The third result connects the convex-set valued maximal operator with matrix $\A_1$ weights.  Given a locally integrably bounded convex-set valued function $F$, we say that it is in convex-set valued $\A_{1,\B}^\K$ if  for some constant $C$ and almost every $x$,
$M_\B F(x) \subset F(x)$.  The infimum of all such constants is denoted by $[F]_{\A_{1,\B}^\K}$.  The following result was proved in~\cite[Corollary~7.4]{bownik-cruz-uribe}. 

\begin{lemma} \label{lemma:convex-A1}
Given a basis $\B$, and a matrix weight $W$, $W\in \A_{1,\B}$ if and only if $W\overline{\BB}\in \A_{1,\B}^\K$.  Further,
$[W]_{\A_{1,\B}} = [W\overline{\BB}]_{\A_{1,\B}^\K}$.
\end{lemma}

\medskip

\begin{proof}[Proof of Theorem~\ref{thm:sharp-off-diag-extrapol-pairs}] 
We will prove this result assuming $p_0<q_0$, so that $1<s<\infty$. The proof when $p_0=q_0$, requires only  minor changes.  Moreover, this case is proved in~\cite[Theorem~9.1]{bownik-cruz-uribe} when $\B=\Q$, the basis of  cubes, and this proof readily adapts to general bases. 
   We divide the proof into four cases: $1<p_0 <q_0<\infty$ and $p<p_0$;  $1<p_0$ and $q_0=\infty$; $1<p_0 < q_0<\infty$ and $p_0<p$; and  $p_0=1$ and $q_0<\infty$. 

\medskip

\textbf{Case I: $\boldsymbol{1<p_0<q_0<\infty}$ and $\boldsymbol{p<p_0}$}. Fix $1<p_0<q_0 <\infty$, and let $p,\,q,\,r,\,s$  be as in the hypotheses.  Let $W$ be a matrix weight such that $W^s \in \A_{r,\B}$.  Define $ P_{W} = N_{W^s} M_\B$; then  $ P_{W} $ is sublinear and monotone, and by Theorem~\ref{thm_6.9},
\[
\| P_{W} \|_{L^r_{\K}(\R^n, W^s)} = \| M_\B \|_{L^r_{\K}(\R^n, W^s)}= K_r([W^s]_{\A_{r,\B}}).
\]
If we define the iteration operator ${R}_W$ by 
\[
\mathcal{R}_W H(x) = \sum_{k=0}^{\infty} 2^{-k} \| P_{W} \|_{L^r_{\K}(\R^n, W^s)}^{-k} P_{W}^k H(x),
\]
then by Theorem~\ref{thm_7.6} it is defined on $ L^r_{\K}(\R^n, W^s) $ and satisfies:
\begin{enumerate} 
    \item[(A$_0$)] $ H(x) \subset \mathcal{R}_W H(x) $;
    \item[(B$_0$)] $ \| \mathcal{R}_W H \|_{L^r_{\K}(\R^n, W^s)} \leq 2 \| H \|_{L^r_{\K}(\R^n, W^s)} $;
    \item[(C$_0$)] $ \mathcal{R}_W H \in \A_1^{\K} $ and 
    $M(\mathcal{R}_W H)(x) \subset  P_{W}(\mathcal{R}_W H)(x) \subset 2K_{r,\B}([W^s]_{\A_{r,\B}})\mathcal{R}_W H(x)$. 
\end{enumerate} 

Now fix $ (f, g) \in \mathcal{F}$; by the definition of $\F$,  $0 < \|f\|_{L^q(\R^n,W)}, \|g\|_{L^p(\R^n,W)} < \infty$.   Define the convex-set valued functions
\[ F(x)=\clconv\{-f(x),f(x)\}, \qquad G(x)=\clconv\{-g(x),g(x)\}.  \]
Then we have that 
\[
 N_WF(x)=|W(x)F(x)|W^{-1}(x)\widebar{\mathbf{B}} = |W(x)f(x)|W^{-1}(x)\widebar{\mathbf{B}}; 
\]
similarly, $N_WG(x)= |W(x)g(x)|W^{-1}(x)\widebar{\mathbf{B}}$.
Define the ellipsoid-valued function
\[
H_1(x) = \bigg( \frac{|W(x) f(x)|}{\| f \|_{L^q(\R^n, W)}} + \bigg[\frac{|W(x) g(x)|}{\| g \|_{L^p(\R^n, W)}}\bigg]^{\frac{p}{q}}
\bigg) W^{-1}(x) \widebar{\mathbf{B}}= h_1(x)W^{-1}(x)\widebar{\mathbf{B}};
\]
then $\| H_1 \|_{L^q(\R^n, W)} = \|h_1\|_{L^q(\R^n)} \leq 2$. 
Define $H_1^s = h_1^s W^{-s} \widebar{\mathbf{B}}$; then
\begin{equation}\label{eq_H1s}
    \| H_1^s \|_{L^r_{\K}(\R^n, W^s)} 
    = \bigg(\int_{\R^n} \left| W^s(x) h^s_1(x)W^{-s}(x)\widebar{\mathbf{B}} \right|^r \,dx \bigg)^{\frac{1}{r}}
    = \bigg( \int_{\R^n} \left| h_1 \right|^q \,dx \bigg)^{\frac{1}{q}s}  \leq 2^s. 
\end{equation}
Thus,  $H_1^s \in L^r_{\K}(\R^n, W^s)$ and $ \mathcal{R}_W H_1^s $ is defined. Since $H_1^s$ is  an ellipsoid-valued function, by Theorem~\ref{thm_7.6}, $\mathcal{R}_W H^s_1$ is also an ellipsoid-valued function. Therefore, there exists a scalar function $ R_W h_1^s $ such that
\[
\mathcal{R}_W H_1^s(x) = \mathcal{R}_W h_1^s W^{-s} \widebar{\mathbf{B}}.
\]
We claim that properties (A$_0$), (B$_0$), and (C$_0$) imply the following:
\begin{enumerate}
    \item [(A$_1$)]  $ h_1(x) \leq \widebar{h}_1(x) $;
    \item [(B$_1$)] $ \| \widebar{h}_1 \|_{L^q(\R^n)} \leq 2^{\frac{1}{s}+1} $;
    \item [(C$_1$)] $ \widebar{h}_1^s W^{-s} \in \A_{1,\B} $ and $[\widebar{h}_1^s W^{-s}]_{\A_{1,\B}} 
    \leq c(d)K_r([W^s]_{\A_{r,\B}}) $.
\end{enumerate}
To prove (A$_1$): it follows  from (A$_0$) that  $ h_1^s(x) \leq R_W h_1^s(x) $. Thus, $h_1(x) \leq \left( R_W h_1^s \right)^{\frac{1}{s}}(x)= \widebar{h}_1(x)$. Inequality (B$_1$) follows from  (B$_0$) and inequality~\eqref{eq_H1s}:
\begin{multline*}
    \| \widebar{h}_1 \|_{L^q(\R^n)} = \left( \int_{\R^n} \left| (\mathcal{R}_W h_1^s)^{\frac{1}{s}}(x) \right|^q \,dx\right)^{\frac{1}{q}} \\
    = \left( \int_{\R^n} \left| W^s(x) \mathcal{R}_W H_1^s(x) \right|^r \,dx \right)^{\frac{1}{rs}} 
    \leq 2^{\frac{1}{s}} \|H_1^s \|_{L^r_{\K}(\R^n,W^s)}^{\frac{1}{s}} \leq 2^{\frac{1}{s}+1}.
\end{multline*}
To prove (C$_1$): from (C$_0$) we have that
    \[
    \mathcal{R}_W H_1^s(x) = \mathcal{R}_W h_1^s W^{-s} \widebar{\mathbf{B}} \in \A_{1,\B}^{\K},
    \]
    and so by Lemma~\ref{lemma:convex-A1},  $\widebar{h}_1^s W^{-s} \in \A_{1,\B}$ and $[\widebar{h}_1^s W^{-s}]_{\A_{1,\B}} 
    \leq c(d) K_r([W^s]_{\A_{r,\B}})$. 

\medskip

We claim that if we define $W_0(x)=\widebar{h}_1(x)^{-\frac{(q_0-q)}{q_0}}W(x)$, then $W_0^s\in \A_{r_0,\B}$ and
\[
[W_0^s]_{\A_{r_0,\B}} \le C(d,r) [W^s]_{\A_{r,\B}}^{\frac{r}{r_0}}  K_r([W^s]_{\A_{r,\B}})^{1-\frac{r}{r_0}}.
\]
To prove $W_0^s\in \A_{r_0,\B}$, note that by (C$_1$), $\widebar{h}_1^sW^{-s}\in \A_{1,\B}$; hence, $\widebar{h}_1^{-s}W^{s}\in \A_{\infty,\B}$. Thus, by Lemma~\ref{matrix-weightsA_r}, $W_0^s=\widebar{h}_1^{-s( \frac{q-q_0}{q_0} )}W^s=\left( \widebar{h}_1^{-s}W^s \right)^{1-\frac{r}{r_0}}\left(W^s \right)^\frac{r}{r_0} \in \A_{r_0,\B}$, and 
\[
[W_0^s]_{\A_{r_0,\B}} \le C(d) [W^s]_{\A_{r,\B}}^{\frac{r}{r_0}} [\widebar{h}_1^{-s}W^s]_{\A_{\infty,\B}}^{1-\frac{r}{r_0}}.
\]
From Lemma~\ref{cor:duality} we have $[\widebar{h}_1^{-s}W^s]_{\A_{\infty,\B}} \approx [\widebar{h}_1^{s}W^{-s}]_{\A_{1,\B}}$, and by (C$_1$),
\[
[W_0^s]_{\A_{r_0,\B}} \lesssim  C(d) [W^s]_{\A_{r,\B}}^\frac{r}{r_0} [\widebar{h}_1^sW^{-s}]_{\A_{1,\B}}^{1-\frac{r}{r_0}} 
\le C(d,r) [W^s]_{\A_{r,\B}}^\frac{r}{r_0}  K_r ( [W^s]_{\A_{r,\B}} )^{1-\frac{r}{r_0}}.
\]

\medskip

We can now prove inequality \eqref{eqn:sharp-extrapol0}. Since $\frac{1}{p_0}-\frac{1}{q_0}=\frac{1}{s'}$ and $p<p_0$,  $q<q_0$. Therefore, by Hölder's inequality with exponent $\frac{q_0}{q}$ and by (B$_1$),
\begin{align*}
    \|f\|_{L^q(\R^n,W)} 
    &= \bigg( \int_{\R^n} |W(x)f(x)|^q \widebar{h}_1(x)^{-\frac{q(q_0-q)}{q_0}} \widebar{h}_1(x)^{\frac{q(q_0-q)}{q_0}}\,dx \bigg)^{\frac{1}{q}} \\
    & \le \bigg( \int_{\R^n} |W(x)f(x)|^{q_0} \widebar{h}_1(x)^{-(q_0-q)} \,dx \bigg)^{\frac{1}{q_0}} 
    \bigg( \int_{\R^n} \widebar{h}_1(x)^q dx \bigg)^{\frac{q_0-q}{qq_0}}\\
    & = \bigg( \int_{\R^n} |W(x)f(x)|^{q_0} \widebar{h}_1(x)^{-(q_0-q)} \,dx \bigg)^{\frac{1}{q_0}} 2^{(\frac{1}{s}+1 )
    ( \frac{q_0-q}{q_0})}.
\end{align*}
By the definition of $h_1$ and (A$_1$), we have that
\[
\frac{|W(x) f(x)|}{\| f \|_{L^q(\R^n, W)}} \leq h_1(x) \leq \widebar{h}_1(x).
\]
Thus,
\[
\int_{\R^n} |W(x)f(x)|^{q_0} \widebar{h}_1(x)^{-(q_0-q)} \,dx \le \|f\|_{L^q(\R^n,W)}^{(q_0-q)} \int_{\R^n} |W(x)f(x)|^q dx = \|f\|_{L^q(\R^n,W)}^{q_0} <\infty.
\]
Again by property (A$_1$) we have 
\[
\bigg[\frac{|W(x) g(x)|}{\| g \|_{L^p(\R^n, W)}}\bigg]^{\frac{p}{q}} \leq h_1(x) \leq \widebar{h}_1(x);
\]
thus,
\[
|W(x)g(x)|^{p_0} \le \widebar{h}_1(x)^{\frac{q}{p}p_0} \|g\|^{p_0}_{L^p(\R^n,W)}.
\]
So, by our assumption \eqref{eqn:sharp-extrapol0}, the definition of $W_0$, (B$_1$), and since $N_{p_0,q_0}$ is increasing, we get
\begin{align*}
    \|f\|_{L^q(\R^n,W)} & \le 2^{( \frac{1}{s}+1 )( \frac{q_0-q}{q_0} )} \left( \int_{\R^n} |W_0(x)f(x)|^{q_0}dx \right)^{\frac{1}{q_0}} \\
    & \le 2^{( \frac{1}{s}+1 )( \frac{q_0-q}{q_0} )} N_{p_0,q_0}\left([W_0^s]_{\A_{r_0,\B}} \right) \left( \int_{\R^n} |W_0(x)g(x)|^{p_0}dx \right)^{\frac{1}{p_0}} \\
    & \le 2^{( \frac{1}{s}+1 )( \frac{q_0-q}{q_0} )} N_{p_0,q_0}\left([W_0^s]_{\A_{r_0,\B}} \right) \|g\|_{L^p(\R^n,W)} \left( \int_{\R^n} \widebar{h}_1(x)^{\frac{q}{p}p_0 - \frac{(q_0-q)}{q_0}p_0} dx \right)^\frac{1}{p_0} \\
    & = 2^{( \frac{1}{s}+1 )( \frac{q_0-q}{q_0} )} N_{p_0,q_0}\left([W_0^s]_{\A_{r_0,\B}} \right) \|g\|_{L^p(\R^n,W)} \|\widebar{h}_1 \|_{L^q(\R^n)}^\frac{q}{p_0}\\
    & \le 2^{( \frac{1}{s}+1 )( \frac{q_0-q}{q_0} )} N_{p_0,q_0}\left([W_0^s]_{\A_{r_0,\B}} \right) \|g\|_{L^p(\R^n,W)} 
    2^{( \frac{1}{s}+1 )\frac{q}{p_0}}\\
    & \le 2^{( \frac{1}{s}+1 )\frac{q}{p} } 
    N_{p_0,q_0}\big(C(d,r) [W^s]_{\A_{r,\B}}^{\frac{r}{r_0}}  K_r([W^s]_{\A_{r,\B}})\big)^{1-\frac{r}{r_0}} \|g\|_{L^p(\R^n,W)};
\end{align*}
in the third inequality we used the fact that $\big( \frac{q}{p}-\frac{(q_0-q)}{q_0} \big)p_0=q$. 

\medskip

\textbf{Case II: $\boldsymbol{1<p_0<\infty}$ and $\boldsymbol{q_0=\infty}$}. Fix $1< p_0 <\infty$ and let $p,q,r$ be as in the hypotheses. Note that in this case, $s=p'_0$ and $r_0=\infty$. Fix a matrix weight $W$ such that $W^{p'_0} \in \A_{r,\B}$. The proof is similar to, but
simpler than, the proof of Case I. 

If we define $H_1$, $h_1$, and $\widebar{h}_1$ as in the previous case, then conditions (A$_1$), (B$_1$), and (C$_1$) hold.
By the definition of $h_1$ and (A$_1$),
\begin{equation} \label{eqn:matrix-for-infty}
    |W(x)f(x)|\widebar{h}_1(x)^{-1} \le |W(x)f(x)|h_1(x)^{-1} \le \|f\|_{L^q(\R^n,W)}<\infty.
\end{equation}
Hence, if we define $W_0(x)=\widebar{h}_1(x)^{-1}W(x)$, we get that
\[
\|f\|_{L^\infty (\R^n,W_0)} \le \|f\|_{L^q(\R^n,W)}<\infty.
\]
By (C$_1$) we have $\widebar{h}_1(x)^{p'_0}W(x)^{-p'_0}\in \A_{1,\B}$; thus, $W_0(x)^{p'_0}=\widebar{h}_1(x)^{-p'_0}W(x)^{p'_0}\in \A_{\infty,\B}$ and
\[
[\widebar{h}_1^{-p'_0}W^{p'_0}]_{\A_{\infty,\B}} \le C(d) K_r \big([W^{p'_0}]_{\A_{r,\B}} \big).
\]

Therefore, by  assumption \eqref{eqn:Sharp-extrapol1} and (B$_1$),
\begin{multline*}
    \|f\|_{L^q(\R^n,W)} 
     = \left(\int_{\R^n} |W_0(x)f(x)|^q \widebar{h}_1(x)^q dx \right)^\frac{1}{q} \\
     \le \|f\|_{L^\infty(\R^n,W_0)} \|\widebar{h}_1\|_{L^q(\R^n)} 
     \le 2^{\frac{1}{p'_0}+1}N_{p_0,\infty} \big( [W_0^{p'_0}]_{\A_{\infty,\B}} \big) \|g\|_{L^{p_0}(\R^n,W_0)}.
\end{multline*}
If we argue as we did to prove \eqref{eqn:matrix-for-infty}, we have that
\[
|\widebar{h}_1(x)^{-1}W(x)g(x)|\le \widebar{h}_1(x)^{\frac{q}{p}-1}\|g\|_{L^p(\R^n,W)}.
\]
Thus, 
\begin{multline*}
  \|g\|_{L^{p_0}(\R^n,W_0)} 
  =  \bigg(\int_{\R^n} |\widebar{h}_1(x)^{-1}W(x)g(x)|^{p_0} dx \bigg)^\frac{1}{p_0} 
  \le \bigg( \int_{\R^n} \|g\|_{L^p(\R^n,W)}^{p_0} \widebar{h}_1(x)^{(\frac{q}{p}-1 )p_0} dx \bigg)^\frac{1}{p_0} \\
    = \|g\|_{L^p(\R^n,W)} \bigg( \int_{\R^n} \widebar{h}_1(x)^q dx \bigg)^\frac{1}{p_0} 
    \le \|g\|_{L^p(\R^n,W)} 2^{(\frac{1}{p'_0}+1 )\frac{q}{p_0}};
\end{multline*}
for the middle equality we use that since $q_0=\infty$, $(\frac{q}{p}-1 )p_0=q$.
If we combine this with the previous inequality and use the fact that $N_{p_0,\infty}$ is increasing, we get
\begin{multline*}
    \|f\|_{L^q(\R^n,W)} 
     \le 2^{(\frac{1}{p'_0}+1 )(1+ \frac{q}{p_0} )} N_{p_0,\infty} 
    \big( [W_0^{p'_0}]_{\A_{\infty,\B}} \big) \|g\|_{L^p(\R^n,W)} \\
     \le 2^{( \frac{1}{p'_0}+1 )\frac{q}{p}} N_{p_0,\infty} 
     \big(C(d) K_r \big( [W^{p'_0}]_{\A_{r,\B}} \big) \big) \|g\|_{L^p(\R^n,W)}.
\end{multline*}

\medskip

\textbf{Case III: $\boldsymbol{1<p_0 < q_0<\infty}$ and $\boldsymbol{p_0<p}$}. Fix $1<p_0<q_0<\infty$, and define $p,\,q,\,r,\,s$  as in the hypotheses. Then $q_0<q$. Fix a matrix weight $W$ such that $W^s \in \A_{r,\B}$. Then, by Corollary~\ref{cor:duality},  $W^{-s}\in \A_{r',\B}$ and so $ M_\B $ is bounded on $ L^{r'}_{\K}(\R^n, W^{-s}) $.  If we define
\[
M'_s H(x) = W^{-s}(x) M_\B(W^s H)(x),
\]
then
\[
\| M'_s H \|_{L^{r'}_{\K}(\R^n)} 
= \| M_\B(W^s H) \|_{L^{r'}_{\K}(\R^n, W^{-s})} 
\leq K_{r'}([W^{-s}]_{\A_{r',\B}}) \| H \|_{L^{r'}_{\K}(\R^n)}.
\]

Let $ P'_{I} = N_{I} M'_s $;  then $ P'_{I} $ is sublinear, monotone and its operator norm is the same as $M'_s$.  Hence, by Theorem \ref{thm_7.6}, if we define the iteration operator
\[
\mathcal{R}'_{I} H(x) = \sum_{k=0}^{\infty} 2^{-k} \| P'_{I} \|_{L^r_{\K}(\R^n)}^{-k} (P'_{I})^k H(x),
\]
then the following hold:
\begin{enumerate}
    \item[(A$'_0$)] $ H(x) \subset \mathcal{R}'_{I} H(x) $;
    \item[(B$'_0$)] $ \| \mathcal{R}'_{I} H \|_{L^{r'}_{\K}(\R^n)} \leq 2 \| H \|_{L^{r'}_{\K}(\R^n)} $;
    \item[(C$'_0$)] $ W^s \mathcal{R}'_{I} H \in \A_1^{\K} $ and $ M(W^s \mathcal{R}'_{I} H)(x) \subset K_{r'}([W^{-s}]_{\A_{r',\B}}) W^s \mathcal{R}'_{I} H(x) $.
\end{enumerate}

Fix $(f,g)\in \mathcal{F}$; again we have that $0<\|f\|_{L^q(w)}, \,\|g\|_{L^p(w)} < \infty$.  Since $\| f \|_{L^q(\R^n, W)}= \| |Wf|^s \|^{\frac{1}{s}}_{L^r(\R^n)} $, by  duality there exists a function $h_2(x)\in L^{r'}(\R^n)$, $\|h_2\|_{L^{r'}(\R^n)}=1$, such that
\[
\|f\|_{L^q(\R^n,W)}=\bigg(\int_{\R^n} |W(x)f(x)|^s h_2(x)dx \bigg)^{\frac{1}{s}}.
\]
Let $H_2(x)=h_2(x)\widebar{\mathbf{B}}$; then $H_2\in L^{r'}_{\K}(\R^n)$, $\|H_2\|_{L^{r'}_{\K}(\R^n)}=1$,. Moreover, since $H_2$ is a ball-valued function, by Theorem~\ref{thm_7.6}, $\mathcal{R}'_{I}H_2$ is also ball-valued. Therefore, there exists a scalar function $\widebar{h}_2= \mathcal{R}'_{I}h_2$ such that
\[
\mathcal{R}'_{I}H_2(x) =\mathcal{R}'_{I}h_2(x)\widebar{\mathbf{B}}.
\]
We now claim that the following hold:
\begin{enumerate}
    \item [(A$_1'$)]  $ h_2(x) \leq \widebar{h}_2(x) $;
    \item [(B$_1'$)] $ \| \widebar{h}_2 \|_{L^{r'}(\R^n)} \leq 2 $;
    \item [(C$_1'$)] $ \widebar{h}_2 W^{s} \in \A_{1,\B} $ and $[\widebar{h}_2 W^{s}]_{\A_{1,\B}} \leq c(d)K_{r'}([W^{-s}]_{\A_{r',\B}}) $.
\end{enumerate}
Inequality (A$_1'$) follows from (A$_0'$) and the definition of $\widebar{h}_2$.  Inequality (B$_1'$) follows from (B$_0'$):
\[ \| \widebar{h}_2 \|_{L^{r'}(\R^n)} = 
\| \mathcal{R}'_{I} H_2 \|_{L^{r'}_{\K}(\R^n)} \leq 2 \| H_2 \|_{L^{r'}_{\K}(\R^n)} = 2. \]
Finally, (C$_1'$) follows from (C$_0'$) and Lemma~\ref{lemma:convex-A1}.

We can now begin to estimate~\eqref{eqn:Sharp-extrapol1}.  By (A$'_1$) and Hölder's inequality with exponents $r_0=\frac{q_0}{s}$ and $r'_0=\frac{p'_0}{s}$,
\begin{align*}
    \|f\|_{L^q(\R^n,W)} &= \left( \int_{\R^n} |W(x)f(x)|^s h_2(x)dx \right)^\frac{1}{s} \\
    &\le \left(\int_{\R^n} |W(x)f(x)|^s \widebar{h}_2(x)^{\frac{(q-q_0)}{q_0(q-s)}s} h_2(x)^{\frac{(q_0-s)q}{q_0(q-s)}} \right)^\frac{1}{s} \\
    &\le \left(\int_{\R^n} |W(x) \widebar{h}_2(x)^{\frac{(q-q_0)}{q_0(q-s)}}f(x)|^{q_0}dx \right)^{\frac{1}{q_0}} \left( \int_{\R^n} h_2(x)^{\frac{q}{q-s}} \right)^\frac{1}{sr'_0} \\
    & = \left(\int_{\R^n} |W(x) \widebar{h}_2(x)^{\frac{(q-q_0)}{q_0(q-s)}}f(x)|^{q_0}dx \right)^{\frac{1}{q_0}};
\end{align*}
the last equality holds since $\frac{q}{q-s}=r'$ and By $\|h_2\|_{L^{r'}(\R^n)}=1$.  This last term is finite:  by Hölder's inequality with $\frac{q}{q_0}$, $( \frac{q}{q_0})'=\frac{q}{q-q_0}$, and by (B$'_1$),
\begin{multline*}
    \|f\|_{L^q(\R^n,W)} 
     \le \left(\int_{\R^n} |W(x)f(x)|^{q_0} \widebar{h}_2(x)^{\frac{(q-q_0)}{q-s}} dx \right)^{\frac{1}{q_0}}  \\
     \le \left( \int_{\R^n} |W(x)f(x)|^q dx \right)^\frac{1}{q} \left( \int_{\R^n} \widebar{h}_2(x)^{r'}dx \right)^\frac{q-q_0}{qq_0} 
     \le 2^{r'\frac{(q-q_0)}{qq_0}} \|f\|_{L^q(\R^n,W)} <\infty.
\end{multline*}

Define $W_0(x)=\widebar{h}_2(x)^{\frac{(q-q_0)}{q_0(q-s)}}W(x)$; we claim that $W_0^s \in \A_{r_0,\B}$ with
\[
[W_0^s]_{\A_{r_0,\B}} \le C(d,r)  K_{r',\B}\big( [W^{-s}]_{\A_{r',\B}} \big)^{1-\frac{r'}{r'_0}} [W^s]_{\A_{r,\B}}^\frac{r'}{r'_0}.
\]
To show this, note that by (C$'_1$), $\widebar{h}_2W^s \in \A_{1,\B}$. Since $r'= \frac{p'}{s}=\frac{q}{q-s}$ and $r'_0=\frac{p'_0}{s}=\frac{q_0}{q_0-s}$, we get $1-\frac{r'}{r'_0}=\frac{s(q-q_0)}{q_0(q-s)}$. Thus, by Lemma~\ref{matrix-weightsA_r}, 
\[ W_0^s(x)=\widebar{h}_2(x)^\frac{s(q-q_0)}{q_0(q-s)} W(x)^s 
= \left( \widebar{h}_2(x) W(x)^s \right)^{1-\frac{r'}{r'_0}} \left( W(x)^s \right)^\frac{r'}{r'_0}\in \A_{r_0,\B}, \]
and 
\[
[W_0^s]_{\A_{r_0,\B}} 
\le C(d) [\widebar{h}_2W^s]_{\A_{1,\B}}^{1-\frac{r'}{r'_0}}[W^s]_{\A_{r,\B}}^\frac{r'}{r'_0} 
\le C(d,r)  K_{r',\B} \big( [W^{-s}]_{\A_{r',\B}} \big)^{1-\frac{r'}{r'_0}} [W^s]_{\A_{r,\B}}^{\frac{r'}{r'_0}}.
\]

\medskip

We can now continue our estimate of~\eqref{eqn:Sharp-extrapol1}.
By~\eqref{eqn:sharp-extrapol0}, Hölder's inequality with exponents $\frac{p}{p_0}$, $( \frac{p}{p_0} )'=\frac{p}{p-p_0}$, 
\begin{align*}
    \|f\|_{L^q(\R^n,W)} 
    & \le \left(\int_{\R^n} |W_0(x) f(x)|^{q_0}dx \right)^\frac{1}{q_0}  \\
    & \le N_{p_0,q_0}([W_0^s]_{\A_{r_0,\B}})\left( \int_{\R^n} |W_0(x)g(x)|^{p_0}dx \right)^\frac{1}{p_0}\\
    & \le N_{p_0,q_0}([W_0^s]_{\A_{r_0,\B}}) \left(\int_{\R^n} |W(x)g(x)|^p dx \right)^\frac{1}{p} \left( \int_{\R^n} \widebar{h}_2(x)^{\frac{p_0(q-q_0)}{q_0(q-s)}\frac{p}{p-p_0}} dx\right)^\frac{p-p_0}{pp_0};\\
    \intertext{since $\frac{p_0(q-q_0)}{q_0(q-s)}\frac{p}{p-p_0}=r'$, by (B$'_1$) and  since $N_{p_0,q_0}$ is increasing,}
    & = N_{p_0,q_0}([W_0^s]_{\A_{r_0,\B}}) \|g\|_{L^p(\R^n,W)}\|\widebar{h}_2 \|_{L^{r'}(\R^n)}^{r'\frac{(p-p_0)}{pp_0}}\\
    & \le 2^{r'\frac{p-p_0}{pp_0}} N_{p_0,q_0}
    \big( C(d,r) K_{r'}\big([W^{-s}]_{\A_{r',\B}} \big)^{1-\frac{r'}{r'_0}} 
    [W^s]_{\A_{r,\B}}^{\frac{r'}{r'_0}} \big) \|g\|_{L^p(\R^n,W)}.
\end{align*}

\medskip

\textbf{Case IV: $\boldsymbol{p_0=1}$ and $\boldsymbol{q_0<\infty}$}. Fix $1<q_0<\infty$ and let $p,q,r$ be as in the hypotheses. Note that in this case $s=q_0$ and $r_0=1$. Fix a matrix weight $W$ such that $W^{q_0}\in \A_{r,\B}$.  The proof is very similar to the proof in Case III.  Again by duality, there exists  $h_2 \in L^{r'(\R^n)}$, $\|h_2\|_{L^{r'}(\R^n)}=1$, such that
\[  \|f\|_{L^q(\R^n,W)} = \left(\int_{\R^n} |W(x)f(x)|^{q_0} h_2(x) dx \right)^\frac{1}{q_0}. \]
Define $H_2$ and $\widebar{h}_2= \mathcal{R}'_Ih_2$ as before. Then (A$'_1$), (B$'_1$), and (C$'_1$) hold. 
Define the matrix $W_0(x) = \widebar{h}_2(x)^\frac{1}{q_0}W(x)$; then by (C$'_1$) $W_0^{q_0} \in \A_{1,\B}$ and
\[
[W_0^{q_0}]_{\A_{r,\B}} \le C(d) K_{r'} \left( [W^{-q_0}]_{\A_{r',\B}} \right).
\]
Moreover, by  the duality relationship, Hölder's inequality with exponent $\frac{q}{q_0}$, and (B$'_1$),
\begin{multline*}
    \|f\|_{L^{q_0}(\R^n,W_0)} 
    = \left(\int_{\R^n} |W(x)f(x)|^{q_0} \widebar{h}_2(x) dx \right)^\frac{1}{q_0}\\
         \le \left( \int_{\R^n} |W(x)f(x)|^q dx \right)^\frac{1}{q} \left(\int_{\R^n} \widebar{h}_2(x)^{r'}dx \right)^\frac{1}{r's} 
    \leq 2^{\frac{1}{s}} \|f\|_{L^q(\R^n,W)} < \infty . 
\end{multline*}

We can now prove~\eqref{eqn:Sharp-extrapol1}.  By the duality relation, (A$_1'$), our assumption \eqref{eqn:sharp-extrapol0} (recalling that $p_0=1$),  Hölder's inequality, and since $N_{1,q_0}$ is increasing, we have
that
\begin{align*}
    \|f\|_{L^q(\R^n,W)} & \le \left( \int_{\R^n} |W_0(x)f(x)|^{q_0} dx \right)^{\frac{1}{q_0}} \\
    & \le N_{1,q_0} \left([W_0^{q_0}]_{\A_{r_0,\B}}\right) \int_{\R^n} |W_0(x)g(x)|dx;\\
    \intertext{by H\"older's inequality with exponent $p$, (B$'_1$), and since $N_{1,q_0}$ is increasing,}
    & = N_{1,q_0} \left([W_0^{q_0}]_{\A_{r_0,\B}}\right) \int_{\R^n} |W(x)g(x)| \widebar{h}_2(x)^{\frac{1}{q_0}}dx \\
    & \le N_{1,q_0} \left([W_0^{q_0}]_{\A_{r_0,\B}}\right) \|g\|_{L^p(\R^n,W)} \left( \int_{\R^n} \widebar{h}_2(x)^{\frac{p'}{q_0}}dx \right)^\frac{1}{p'} \\
    & = N_{1,q_0} \left([W_0^{q_0}]_{\A_{r_0,\B}}\right) \|g\|_{L^p(\R^n,W)} \|\widebar{h}_2 \|_{L^{r'}(\R^n)}^\frac{1}{q_0} \\
    & \le 2^\frac{1}{q_0} N_{1,q_0} \big( C(d) K_{r'} \big( [W^{-q_0}]_{\A_{r',\B}} \big) \big)  \|g\|_{L^p(\R^n,W)}.
\end{align*}
\end{proof}

\medskip

\begin{proof}[Proof of Theorem~\ref{thm:off-diag-extrapol-sharp-cubes}]
To prove this result, we first estimate the constant in Theorem~\ref{thm:sharp-off-diag-extrapol-pairs} when $\B$ is the basis of cubes $\Q$.  
In this case, by Theorem \ref{thm_6.9} we have that for $W^s\in\A_{r}$,
\[
\|M\|_{L^r_{\K}(\R^n,W^s)} \le C(n,d,r)[W^s]^{r'}_{\A_r}.
\]
Hence, if $p<p_0$, since $N_{p_0,q_0}$ is increasing,
\begin{align*} 
    N_{p,q}(\B,p_0,q_0, [W^s]_{\A_{r,\B}}) 
    & = 2^{( \frac{1}{s}+1 )\frac{q}{p} } N_{p_0,q_0}\big(C(d,r) [W^s]_{\A_{r,\B}}^{\frac{r}{r_0}}  
    K_r( [W^s]_{\A_{r,\B}})^{1-\frac{r}{r_0}}\big) \\
    & \leq 2^{( \frac{1}{s}+1 )\frac{q}{p} } N_{p_0,q_0}\big(C(n,d,p,q,p_0,q_0)
    [W^s]_{\A_{r,\B}}^{\frac{r}{r_0}+(1-\frac{r}{r_0})r'}\big) \\
    & = 2^{( \frac{1}{s}+1 )\frac{q}{p} } N_{p_0,q_0}\big(C(n,d,p,q,p_0,q_0)
    [W^s]_{\A_{r,\B}}^{\frac{r'}{r_0'}}\big).
\end{align*}
On the other hand, if $p>p_0$, by Corollary~\ref{cor:duality},  $W^{-s}\in \A_{r'}$, and so
\[
\|M\|_{L^{r'}_{\K}(\R^n,W^{-s})} \le C(n,d,p)[W^{-s}]^{r}_{\A_{r'}}.
\]
Then, by a very similar argument, we have that
\begin{multline*}
    N_{p,q}(\B,p_0,q_0, [W^s]_{\A_{r,\B}}) 
    = 2^{r'\frac{p-p_0}{pp_0}} N_{p_0,q_0}\big( C(d,r) [W^s]_{\A_{r,\B}}^{\frac{r'}{r'_0}}  
    K_{r'}([W^{-s}]_{\A_{r',\B}}) )^{1-\frac{r'}{r'_0}} \big) \\
    \leq 2^{r'\frac{p-p_0}{pp_0}}N_{p_0,q_0}\big(C(n,d,p,q,p_0,q_0)
    [W^s]_{\A_{r,\B}}^{\frac{r}{r_0}}\big).
\end{multline*}
If we combine these two estimates, we get the constant that appears in Theorem~\ref{thm:off-diag-extrapol-sharp-cubes}.

To complete the proof we need to construct an appropriate family of extrapolation pairs.  This was done in~\cite[Theorem~1.5]{bownik-cruz-uribe} and we sketch the changes to apply to the off-diagonal setting here.  Suppose $T$ is a linear operator such that~\eqref{eqn:odec1} holds for some $p_0$ and every $W_0$ such that $W_0^s\in \A_{r_0}$.  Then, by letting $W_0$ be the diagonal matrix $w_0I_d$, where $w_0^s$ is in scalar $\A_{r_0}$, we see that $T : L^{p_0}(w_0)\to L^{q_0}(w_0)$.  Hence, by Theorem~\ref{thm:sharp_constant-scalar} below, $T: L^{p}(w)\to L^{q}(w)$, for any $w$ such that $w^s \in \A_r$.  

For all $1\leq p<\infty$ and any matrix weight $W$, by~\cite[Proposition~3.6]{MR3544941}, $L_c^\infty(\R^n,\R^d)$ is dense in $L^p(\R^n,W)$. By~\cite[Corollary~6.7]{bownik-cruz-uribe}, if $W^s\in \A_r$, then $|W^s|_{\op}= |W|_{\op}^s$ is a scalar weight in $\A_r$.  Therefore, we have that for any $f\in L_c^\infty(\R^n,\R^d)$,
\[ \|Tf\|_{L^q(\R^n,W)} \leq \|Tf\|_{L^q(\R^n,|W|_{\op})} \lesssim \|f_{\R^n, L^q(|W|_{\op})}< \infty. \]
We can now apply Theorem~\ref{thm:sharp-off-diag-extrapol-pairs} to the family $\F = \{ (|Tf|,|f|) : f \in L_c^\infty(\R^n,\R^d)\}$ to prove 
that inequality \eqref{eqn:odec2} holds for all $f\in L_c^\infty(\R^n,\R^d)$.  To prove that it holds for all $f\in L^p(\R^n,W)$, we use a standard approximation argument.
\end{proof}

\appendix 

\section{Sharp constant scalar extrapolation}
\label{section:scalar}

In this appendix we prove the scalar version of Theorem~\ref{thm:sharp-off-diag-extrapol-pairs}.  This result was originally proved in~\cite{MR2652182}; here we give a simpler proof based on the proof for the sharp constant extrapolation theorem in the diagonal case in~\cite[Section~3.4]{MR2797562}.  As we noted in the Introduction, the proof of Theorem~\ref{thm:sharp-off-diag-extrapol-pairs} is modeled on this scalar result and its proof can serve as a guide to the larger ideas in the matrix case.  After the proof, we make a few comments about the technical obstacles involved in generalizing from the scalar to the matrix case.

We first give some definitions.  Let $\B$ be a basis; define the scalar maximal operator $M_\B$ by 
\[
M_\B f(x)=\sup_{B\in \B} \avgint_B |f(y)| \, dy \cdot \chi_B(x).
\]
Given a basis $\B$, A scalar weight is said to satisfy $w\in \A_{p,\B}, 1 \leq p \leq \infty$, if
\[
[w]_{\A_{p,\B}}=\sup_{B\in \B} |B|^{-1}\|w\chi_B\|_{L^p(\R^n)} \|w^{-1}\chi_B\|_{L^{p'}(\R^n)}<\infty.
\]
Note that it is immediate from the definition that $[w]_{\A_{p,\B}}=[w^{-1}]_{\A_{p',\B}}$.
A basis $\B$ is said to be a Muckenhoupt basis if for every $1<p<\infty$ and every $w\in \A_{p,\B}$, $M_\B : L^p(\R^n, w) \to L^p(\R^n,w)$.   As before, we  denote the best constant in the norm inequality by $K_{p,\B}([w]_{\A_{p,\B}})$.  

For $1 \le p \le q \le \infty$ we define the class of scalar weights 
 $\A_{p,q,\B}$ by 
\[
[w]_{\A_{p,q,\B}} = \sup_B |B|^{-1} \|w\chi_B\|_{L^q(\R^n)}\|w^{-1}\chi_B\|_{L^{p'}(\R^n)} <\infty.
\]
An immediate calculation shows that if $r=\frac{q}{s}$, where we define $s$ by  $\frac{1}{p}-\frac{1}{q}=\frac{1}{s'}$, then $r'=\frac{p'}{s}$.  Hence, $w\in \A_{p,q,\B}$ if and only if $w^s\in \A_{r,\B}$.  In fact, we have that $[w]_{\A_{p,q,\B}} = [w^s]_{\A_{r,\B}}^{\frac{1}{s}}$.  When $1<p\leq q<\infty$, this follows immediately from the definitions, since for all $B\in \B$,
\[
\bigg( \avgint_{B} w^q \, dx \bigg)^{\frac{1}{q}} 
\left( \avgint_{B}  w^{-p'} \, dx \right)^{\frac{1}{p'}} = \left( \avgint_{B}  (w^s)^r \, dx \right)^{\frac{1}{r s}} 
\left( \avgint_{B} (w^s)^{-r'} \, dx \right)^{\frac{1}{r' s}}.
\]
A similar argument holds when $p=1$ or $q=\infty$.  Because of this, we will state our result for $\A_{p,q,\B}$ weights, and use this equivalent condition in the proof.

\begin{theorem}\label{thm:sharp_constant-scalar}
Let $\B$ be a Muckenhoupt basis and let $\F$ be a family of extrapolation pairs. Suppose that for some $ p_0, q_0, 1 \leq p_0 \leq q_0 \leq \infty $,  there exists a positive increasing function $ N_{p_0,q_0}$ such that for every weight $w \in A_{p_0,q_0,\B}$,
\begin{equation} \label{eqn:extrapol0}
   \|f\|_{L^{q_0}(\R^n,w_0)}
\leq N_{p_0, q_0} \left( [w_0]_{A_{p_0,q_0,\B}}\right)
\|g\|_{L^{p_0}(\R^n,w_0)}, \quad (f,g) \in \F.
\end{equation}
Then for all $ p,\, q $ such that $ 1 < p \leq q < \infty $, $ \frac{1}{p} - \frac{1}{q} = \frac{1}{p_0} - \frac{1}{q_0}$, and for all weights $w \in \A_{p,q,\B}$,
\begin{equation} \label{eqn:extrapol1}
\|f\|_{L^{q}(\R^n,w)}
\leq N_{p, q}(\B,p_0,q_0,[w]_{\A_{p,q,\B}})
\|g\|_{L^{p}(\R^n,w)}, \qquad (f,g) \in \F,
\end{equation}
where if $q<q_0$,
\begin{equation} \label{eqn:first-const}
N_{p, q}(\B,p_0,q_0,[w]_{\A_{p,q,\B}})
= 2^{( \frac{1}{s}+1 ) \frac{q}{p}} N_{p_0,q_0} \big( ( 2\|M_\B \|_{L^r(\R^n,w^s)})^{\frac{q_0-q}{q_0s}} [w]^{\frac{q}{q_0}}_{\A_{p,q,\B}}  \big),
\end{equation}
and if $q_0 < q $,
\begin{equation} \label{eqn:second-const}
N_{p, q}(\B,p_0,q_0,[w]_{\A_{p,q,\B}})
= 2^{ \frac{q-q_0}{q_0(q-s)} }  N_{p_0,q_0} \big( ( 2\|M_\B \|_{L^{r'}(w^{-s})})^{\frac{q-q_0}{q_0(q-s)}} [w]^{\frac{p'}{p'_0}}_{\A_{p,q,\B}}  \big).
\end{equation}
Here,  $1\leq s<\infty$ is defined by  $\frac{1}{p} - \frac{1}{q} = \frac{1}{p_0}-\frac{1}{q_0}=\frac{1}{s'}$,
\end{theorem}

\begin{remark}
    In the statement of Theorem~\ref{thm:sharp_constant-scalar} we exclude the case $s=\infty$, which corresponds to $\frac{1}{p_0}-\frac{1}{q_0}=1$, as the only values of $p_0$ and $q_0$ for which this can hold are $p_0=1$ and $q_0=\infty$.  For in this case there are no values of $p$ and $q$ to extrapolate to.
\end{remark}

\begin{corollary} \label{cor:scalar-extrapol}
    If $\B=\Q$, the basis of cubes, then the constant in Theorem~\ref{thm:sharp_constant-scalar} becomes
    \[
N_{p, q}(\Q,p_0,q_0,[w]_{\A_{p,q}})
= C(p,q,p_0,q_0)  N_{p_0,q_0} \big( C(n,p,q,p_0,q_0) [w]^{\max \{ \frac{q}{q_0},\frac{p'}{p'_0} \}}_{\A_{p,q}}  \big).
\]
\end{corollary}

\begin{proof}[Proof of Corollary~\ref{cor:scalar-extrapol}]
 when $\B$ is the basis of cubes $\Q$ and $M_\B$ is the Hardy-Littlewood maximal operator, Buckley \cite{buckley1993} proved that 
\[
\|M\|_{L^r(\R^n,w)} \le C(n,r) [w]^{r'}_{\A_r}.
\]
Thus,
\[ \|M\|_{L^r(\R^n,w^s)} \le C(n,r) [w^s]^{r'}_{\A_r}=C(n,r) [w]_{\A_{p,q}}^{p'}, \]
and if we insert this into~\eqref{eqn:first-const}, we see that the exponent on $[w]_{\A_{p,q}}$ is
\[ p'\frac{q_0-q}{q_0s} + \frac{sq}{sq_0} = \frac{p'}{p_0'}.  \]
Similarly, we have that 
\[ 
\|M\|_{L^{r'}(w^{-s})} \le C(n,r) [w^{-s}]^{r}_{\A_{r'}}=C(n,r) [w^s]^r_{\A_r} = C(n,r) [w]^q_{\A_{p,q}},
\]
and inserting this into~\eqref{eqn:second-const} the exponent on $[w]_{\A_{p,q}}$ is $\frac{q}{q_0}$. 
\end{proof}

\medskip

\begin{proof}[Proof of Theorem~\ref{thm:sharp_constant-scalar}]
  The case $p_0 = q_0$ is proved in \cite[Theorem 3.22]{MR2797562}, so we will  assume that $p_0 < q_0$, so that $1<s<\infty$. There are four cases: $1<p_0 <q_0<\infty$ and $p<p_0$;  $1<p_0$ and $q_0=\infty$; $1<p_0 < q_0<\infty$ and $p_0<p$; and  $p_0=1$ and $q_0<\infty$. 
  
\medskip
    
\textbf{Case I: $\boldsymbol{1<p_0<q_0<\infty}$ and $\boldsymbol{p<p_0}$}.
    Fix $1<p_0<q_0 <\infty$, and let $p,\,q,\,r,\,s$  be as in the hypotheses. 
     Fix $w\in A_{p,q,\B}$; then $w^s\in A_{r,\B}$, and so $M_\B$ is bounded on $L^r(\R^n,w^s)$. Therefore, given non-negative $h_1 \in L^r(\R^n,w^s)$ we can define the iteration operator $\Rdf$ by
    \[
    \mathcal{R} h_1 (x)=\sum_{k=0}^{\infty} \frac{M_\B^k h_1(x)}{2^k \|M_\B\|^k_{L^r(\R^n,w^s)}}.
    \]
    It is straightforward to show (see~\cite[Theorem~3.9]{MR2797562}) that $\mathcal{R}h_1$ satisfies the following:
    \begin{enumerate}
        \item [(A$_0$)] $h_1(x) \le \mathcal{R}h_1(x)$;
        \item [(B$_0$)] $\|\mathcal{R}h_1\|_{L^r(\R^n,w^s)} \le 2\|h_1\|_{L^r(\R^n,w^s)}$;
        \item [(C$_0$)] $[\mathcal{R}h_1]_{\A_{1,\B}} \le 2\|M_{\B} \|_{L^r(\R^n,w^s)}$.
    \end{enumerate}
    Fix $(f,g) \in \F$; by the definition of $\F$ (in the scalar case), $0<\|f\|_{L^q(w)}, \,\|g\|_{L^p(w)} < \infty$. 
    Define  
    \[
    h_1(x) = \frac{f(x)}{\|f\|_{L^q(\R^n,w)}} + \frac{g(x)^{p/q} w(x)^{p/q-1}}{\|g\|_{L^p(\R^n,w)}^{p/q}}.
    \]
    Then
    \begin{multline} \label{eqn:h1-est}
        \| h_1 \|_{L^q(\R^n,w)} = \left( \int_{\mathbb{R}^n} (w(x) h_1(x))^q \, dx \right)^{\frac{1}{q}}\\
        \leq \left( \int_{\mathbb{R}^n} \frac{(w(x) f(x))^q}{\|f\|^q_{L^q(\R^n,w)}} dx \right)^{\frac{1}{q}}
        + \left( \int_{\mathbb{R}^n} \frac{(w(x) g(x))^p}{\|g\|^p_{L^p(\R^n,w)}}  dx \right)^{\frac{1}{q}} 
         \leq 2.
    \end{multline}
    Thus, $h_1\in L^q(\R^n,w)$. and by rescaling, $\|h_1\|_{L^q(\R^n,w)}=\|h_1^s\|_{L^r(\R^n,w^s)}^{1/s}$, so $h_1^s \in L^r(\R^n,w^s)$. Define $H_1=(\mathcal{R}h_1^s)^{1/s}$. We claim that $H_1$ satisfies
    \begin{enumerate}
        \item [(A$_1$)] $h_1(x) \leq H_1(x)$;
        \item [(B$_1$)] $\|H_1\|_{L^q(\R^n,w)} \le 2^{\frac{1}{s} + 1}$;
        \item [(C$_1$)] $[H_1^s]_{\A_{1,\B}} \leq 2 \|M_\B\|_{L^r(\R^n,w^s)}$.
    \end{enumerate}
    Property (A$_1$) follows at once from  (A$_0$). Inequality (B$_1$) follows from (B$_0$) and \eqref{eqn:h1-est}. Inequality (C$_1$) follows from (C$_0$). 

    \medskip

 We now claim that $w_0 \in A_{p_0, q_0, \B}$.  We prove this using the definition.  Given \( B \in \B \), by (C$_1$), the left-hand integral in the $A_{p_0, q_0, \B}$ condition becomes
    \begin{equation*}
        \left( \avgint_{B} w_0^{q_0} \, dx \right)^{\frac{1}{q_0}} = \left( \avgint_{B} H_1^{-(q_0 - q)} w(x)^q \, dx \right)^{\frac{1}{q_0}} 
        \leq [H_1^s]_{A_{1,\B}}^{\frac{q_0 - q}{q_0 s}} \left( \avgint_{B} H_1^s \, dx \right)^{-\frac{(q_0-q)}{q_0 s}} \left( \avgint_{B} w(x)^q \, dx \right)^{\frac{1}{q_0}} .
    \end{equation*}
   To estimate the  right-hand integral,  we will use Hölder's inequality for $t=\frac{q_0 p'}{q p'_0}$; $t>1$ since $\frac{q_0}{q}>1$ and $\frac{p'}{p'_0}>1$ since $p<p_0$. Also, since  $\frac{1}{p}-\frac{1}{q} = \frac{1}{p_0}-\frac{1}{q_0}=\frac{1}{s'}$,  $p'_0=\frac{q_0s}{q_0-s}$ and $p'=\frac{qs}{q-s}$. Then
    \[
    t'=\frac{t}{t-1}=\frac{q_0p'}{q_0p'-qp'_0}=\frac{\frac{q_0qs}{q-s}}{\frac{q_0qs}{q-s}-\frac{qq_0s}{q_0-s}}=\frac{\frac{1}{q-s}}{\frac{1}{q-s}-\frac{1}{q_0-s}}=\frac{q_0-s}{q_0-q}.
    \]
    Hence, 
    \begin{align*}
        \avgint_B w_0(x)^{-p'_0} dx 
        & = \avgint_B H_1(x)^{(q_0-q)\frac{p'_0}{q_0}}w(x)^{-q\frac{p'_0}{q_0}} dx \\
        & \le \left( \avgint_{B}  H_1(x)^{\frac{q_0 - q}{q_0} \frac{q_0 s}{q_0 - s} \frac{q_0 - s}{q_0 - q}} \, dx \right)^{\frac{q_0-q}{q_0-s}} \left( \avgint_{B}  w(x)^{-p'} \, dx \right)^{\frac{qp_0'}{q_0 p'}}\\
        & = \left( \avgint_{B} H_1(x)^s \, dx \right)^{\frac{q_0 - q}{q_0 - s}} \left( \avgint_{B} w(x)^{-p'} \, dx \right)^{\frac{qp_0'}{q_0 p'}};
    \end{align*}
    if we take the $1/p_0'$ power, we get
    \[
    \left( \avgint_{B}  w_0^{-p'_0} \, dx \right)^{\frac{1}{p'_0}} \leq \left( \avgint_{B} H_1(x)^s \, dx \right)^{\frac{q_0 - q}{q_0 s}} \left( \avgint_{B}  w(x)^{-p'} \, dx \right)^{\frac{q}{p'q_0}}.
    \]
    Thus,
    \begin{multline*}
        \left( \avgint_{B} w_0(x)^{q_0} \, dx \right)^{\frac{1}{q_0}} \left( \avgint_{B}  w_0^{-p'_0} \, dx \right)^{\frac{1}{p'_0}} \leq [H_1^s]_{A_{1,\B}}^{\frac{q_0 - q}{q_0 s}} \left( \avgint_{B}  w(x)^q \, dx \right)^{\frac{q}{qq_0}} \left( \avgint_{B} w(x)^{-p'} \, dx \right)^{\frac{q}{p'q_0}} \\
        \leq [H_1^s]_{A_{1,\B}}^{\frac{q_0 - q}{q_0 s}} [w]^{\frac{q}{q_0}}_{A_{p,q,\B}} \le \left(2 \|M_\B \|_{L^r(\R^n,w^s)} \right)^\frac{q_0-q}{q_0s}[w]^{\frac{q}{q_0}}_{A_{p,q,\B}} < \infty.
    \end{multline*}
    Since this bound holds for all  $B$, $w_0 \in A_{p_0, q_0, \B}$.
    
    We can now prove~\eqref{eqn:extrapol1}.  By Hölder's inequality with exponent $q_0/q$ and by (B$_1$),
    \begin{align*}
        \|f\|_{L^q(\R^n,w)} & = \left( \int_{\mathbb{R}^n} (w(x) f(x))^q \, dx \right)^{\frac{1}{q}} \\
        & = \left( \int_{\mathbb{R}^n} \left(w(x) H_1(x)^{\frac{q_0 - q}{q_0}} H_1(x)^{-\frac{(q_0 - q)}{q_0}} f(x) \right)^q \, dx \right)^{\frac{1}{q}} \\
        & \leq \left( \int_{\mathbb{R}^n} f(x)^{q_0} H_1(x)^{-\left( q_0 - q \right)} w(x)^q \, dx \right)^{\frac{1}{q_0}} \left( \int_{\mathbb{R}^n} H_1(x)^q w(x)^q \, dx \right)^{\frac{q_0 - q}{qq_0}} \\
        & \le  2^{( \frac{1}{s} + 1 ) ( \frac{q_0 - q}{q_0} )} 
        \left( \int_{\mathbb{R}^n} f(x)^{q_0} H_1(x)^{-(q_0 - q)} w(x)^q \, dx \right)^{\frac{1}{q_0}}.
    \end{align*}

    By (A$_1$), $\frac{f(x)}{\| f \|_{L^q(\R^n,w)}} \leq h_1(x) \leq H_1(x)$; hence,
    \[
    \int_{\mathbb{R}^n} f(x)^{q_0} H_1(x)^{-(q_0 - q)} w(x)^q \, dx \leq \| f \|_{L^q(\R^n,w)}^{q_0-q} \int_{\mathbb{R}^n} f(x)^q w(x)^q  \, dx = \| f \|^{q_0}_{L^q(\R^n,w)} < \infty.
    \]

    Therefore, by \eqref{eqn:extrapol1}, property (B$_1$), and since \( N_{p_0, q_0} \) is increasing, 
    \begin{align*}
        \| f \|_{L^q(\R^n,w)} 
        & \leq 2^{( \frac{1}{s} + 1 ) ( \frac{q_0 - q}{q_0} )} N_{p_0, q_0} \left( [w_0]_{A_{p_0, q_0, \B}} \right) 
        \left( \int_{\mathbb{R}^n} g(x)^{p_0} \left( H_1(x)^{-(q_0 - q)} w(x)^q \right)^{\frac{p_0}{q_0}} \, dx \right)^{\frac{1}{p_0}} \\
        & \leq 2^{( \frac{1}{s} + 1 ) ( \frac{q_0 - q}{q_0} )} N_{p_0, q_0} \left( [w_0]_{A_{p_0, q_0, \B}} \right) \| g \|_{L^p(\R^n,w)} \\
        & \qquad \qquad \times \left( \int_{\mathbb{R}^n} H_1(x)^{\frac{q}{p}p_0 - \frac{q_0 - q}{q_0}p_0} w(x)^{\left(\frac{q}{p}-1\right)p_0 + \frac{q}{q_0}p_0} \, dx \right)^{\frac{1}{p_0}} \\
        & = 2^{( \frac{1}{s} + 1 ) ( \frac{q_0 - q}{q_0} )} N_{p_0, q_0} \left( [w_0]_{A_{p_0, q_0, \B}} \right) \| g \|_{L^p(\R^n,w)} \left( \int_{\mathbb{R}^n} H_1(x)^q  w(x)^q \, dx \right)^{\frac{1}{p_0}} \\
        & \leq 2^{( \frac{1}{s} + 1 ) ( \frac{q_0 - q}{q_0} )} N_{p_0, q_0} 
        \big( (2\|M_\B \|_{L^r(\R^n,w^s)})^{\frac{q_0-q}{q_0s}}[w]_{A_{p, q, \B}}^{\frac{q}{q_0}} \big) 
        \| g \|_{L^p(\R^n,w)} \|H_1\|_{L^q(\R^n,w)}^{\frac{q}{p_0}} \\
        & \leq 2^{(( \frac{1}{s} + 1 ) ( \frac{q_0 - q}{q_0} ) + (\frac{1}{s}+1)\frac{q}{p_0})} N_{p_0, q_0} \big( (2\|M_\B \|_{L^r(\R^n,w^s)})^{\frac{q_0-q}{q_0s}}[w]_{A_{p, q, \B}}^{\frac{q}{q_0}} \big) \| g \|_{L^p(\R^n,w)} \\
        & = 2^{( \frac{1}{s} + 1 ) \frac{q}{p} } 
        N_{p_0, q_0} \big( (2\|M_\B \|_{L^r(\R^n,w^s)})^{\frac{q_0-q}{q_0s}}[w]_{A_{p, q, \B}}^{\frac{q}{q_0}} \big) 
        \| g \|_{L^p(\R^n,w)}.
    \end{align*}

\medskip

\textbf{Case II: $\boldsymbol{1<p_0<\infty}$ and $\boldsymbol{q_0=\infty}$}.
Fix $p_0$ and let $p,\, q,\, r,\, s$ be the same as in the hypotheses.  In this case,  $s=p'_0$ and $r_0=\infty$. Fix a weight $w\in \A_{p,q,\B}$. Define $h_1$ and $H_1$ as in Case I. Let $w_0=H_1^{-1}$. For $1<p_0<\infty$ we have that $[H_1^{p'_0}]_{\A_{1,\B}}=[H_1^{-p'_0}]_{\A_{\infty,\B}}$, and 
\begin{multline*}
    [H_1^{-1}]_{\A_{p_0,\infty,\B}}= \sup_{B\in \B} \esssup_{x\in B} H_1^{-1}(x) \left( \avgint_B H_1^{p'_0}(x)dx \right)^{\frac{1}{p'_0}} \\
    = \sup_{B\in \B} \left(\esssup_{x\in B}H_1^{-p'_0}(x) \right)^{\frac{1}{p'_0}} \left( \avgint_B (H_1^{-p'_0})^{-1}(x) dx \right)^{\frac{1}{p'_0}} =[H_1^{-p'_0}]_{\A_{\infty,\B}}^{\frac{1}{p'_0}}.
\end{multline*}
Thus, by (C$_1$),
\[
[H_1^{-1}]_{\A_{p_0,\infty,\B}}=[H_1^{-p'_0}]_{\A_{\infty,\B}}^{\frac{1}{p'_0}} = [H_1^{p'_0}]_{\A_{1,\B}}^{\frac{1}{p'_0}} \le \left( 2\|M_\B \|_{L^r(\R^n,w^s)} \right)^\frac{1}{p'_0} < \infty.
\]
By the definition of $h_1$ and (A$_1$),
\[
f(x)H_1^{-1}(x) \le \|f\|_{L^q(\R^n,w)}<\infty.
\]
Then,
\[
\|f\|_{L^{\infty}(\R^n,w_0)} \le \|f\|_{L^q(\R^n,w)} <\infty.
\]
Therefore, by our assumption \eqref{eqn:extrapol0},
\begin{align*}
    \|f\|_{L^q(\R^n,w)} 
     & =\left( \int_{\mathbb{R}^n} \left(w(x) H_1(x) H_1^{-1}(x) f(x) \right)^q dx \right)^{\frac{1}{q}} \\
     & \le \|f\|_{L^{\infty}(\R^n,w_0)} \|H_1\|_{L^q(\R^n,w)} \\
     & \le 2^{\frac{1}{p'_0}+1} N_{p_0,\infty} \big( [w_0]_{\A_{p_0,\infty,\B}} \big) \|g\|_{L^{p_0}(w_0)} \\
     &\le 2^{\frac{1}{p'_0}+1} N_{p_0,\infty} \big( [w_0]_{\A_{p_0,\infty,\B}} \big) \left( \int_{\mathbb{R}^n} 
     \left( H_1^{-1}(x)g(x) \right)^{p_0} dx \right)^{\frac{1}{p_0}}.
\end{align*}
Again by the definition of $h_1$ and (A$_1$), we have that
\[
\frac{g(x)^{p/q}w(x)^{p/q-1}}{\|g\|_{L^p(\R^n,w)}^{p/q}} \le h_1(x) \le H_1(x),
\]
which in turn implies that
\[
g(x) \le \|g\|_{L^p(\R^n,w)} H_1(x)^{\frac{q}{p}}w(x)^{\frac{q}{p}-1}.
\]
Hence, 
\begin{multline*}
    \left(\int_{\mathbb{R}^n} \left( H_1(x)^{-1} g(x) \right)^{p_0} dx \right)^{\frac{1}{p_0}} 
     \le \|g\|_{L^p(\R^n,w)} \left(\int_{\mathbb{R}^n} \left(H_1(x)w(x) \right)^{( \frac{q}{p}-1 )p_0} dx \right)^{\frac{1}{p_0}} \\
     \le \|g\|_{L^p(\R^n,w)} \|H_1 \|_{L^q(\R^n,w)}^{\frac{q}{p_0}}
     \le 2^{( \frac{1}{p'_0}+1 )\frac{q}{p_0}} \|g\|_{L^p(\R^n,w)}.
\end{multline*}
If we combine these two inequalities, we get
\[
\|f\|_{L^q(\R^n,w)} \le 2^{( \frac{1}{p'_0}+1 )\frac{q}{p} } N_{p_0,\infty}\big( 
( 2\|M_\B \|_{L^r(w^{p'_0})} )^{\frac{1}{p'_0}} \big) \|g\|_{L^p(\R^n,w)}.
\]

\medskip

\textbf{Case III: $\boldsymbol{1<p_0 < q_0<\infty}$ and $\boldsymbol{p_0<p}$}. Fix $1< q_0 <\infty$, and let $p,\,q\,r,\,s$ be as in the hypotheses.  Fix  $w\in A_{p,q,\B}$. Since $w^s \in \A_{r,\B}$, $w^{-s}\in \A_{r',\B}$.  Therefore,  $M_\B$ is  bounded on $L^{r'}(w^{-s})$, so we can define the  iteration operator
\[ 
\mathcal{R}' h_2(x) = \sum_{k=1}^{\infty} \frac{(M_\B)^k h_2(x)}{2^k\|M_\B\|^k_{L^{r'}(w^{-s})}}.
\]
As before, we have  that $\Rdf'$ has the following properties:
\begin{enumerate}
    \item [(A$_0'$)] $h_2(x)\le \mathcal{R}'h_2(x)$;
    \item [(B$_0'$)] $\|\mathcal{R}'h_2\|_{L^{r'}(w^{-s})} \le 2 \|h_2\|_{L^{r'}(w^{-s})}$;
    \item [(C$_0'$)] $[\mathcal{R}'h_2]_{A_{1,\B}} \le 2 \|M_\B \|_{L^{r'}(w^{-s})}$.
\end{enumerate}

Fix $(f,g)\in \mathcal{F}$; again we have that $0<\|f\|_{L^q(w)}, \,\|g\|_{L^p(w)} < \infty$.  Since $f^s \in L^r(\mathbb{R}^n,w^s)$, by  duality there exists $h_2\in L^{r'}(\mathbb{R}^n)$, $\|h_2\|_{L^{r'}(\mathbb{R}^n)}=1$, such that
\begin{equation} \label{eqn:dual-caseIII}
\| f \|_{L^q(\R^n,w)} = \| f^s \|_{L^r(\R^n,w^s)}^{\frac{1}{s}} = \left( \int_{\mathbb{R}^n} (w(x)^s f(x)^s) h_2(x) \, dx \right)^{\frac{1}{s}}.
\end{equation}
Since
\[
\| h_2 w^{s} \|_{L^{r'}(w^{-s})} = \| h_2 \|_{L^{r'}(\mathbb{R}^n)}=1,
\]
we can define  $H_2 = \mathcal{R}'(h_2 w^{s}) w^{-s}$, which has the following properties:
\begin{enumerate}
    \item [(A$_1'$)] $h_2(x) \leq H_2(x)$;
    \item [(B$_1'$)] $\| H_2 \|_{L^{r'}(\mathbb{R}^n)} \leq 2$;
    \item [(C$_1'$)] $[H_2w^s]_{\A_{1,\B}} \le 2 \|M_\B \|_{L^{r'}(w^{-s})}$.
\end{enumerate}
Inequality (A$_1'$) follows from  the definition of $H_2$ and (A$_0'$). Inequality (B$_1'$), follows from (B$_0'$) since 
\[
\| H_2 \|_{L^{r'}(\mathbb{R}^n)} = \| \mathcal{R}'(h_2 w^{s}) \|_{L^{r'}(w^{-s})} \le 2 \| h_2 w^{s} \|_{L^{r'}(w^{-s})} = 2.
\]
Finally, (C$_1'$) is a consequence of (C$_0'$).

\medskip

We claim that $w_0(x) = w(x) H_2(x)^{\frac{q_0 - q}{q_0 (q - s)}} \in A_{p_0, q_0, \B}$ and
\[
[w_0]_{A_{p_0,q_0,\B}} \le \big( 2\|M_\B \|_{L^{r'}(w^{-s})}\big)^{\frac{q-q_0}{q_0(q-s)}} [w]_{\A_{p,q,\B}}^{\frac{p'}{p'_0}}.
\]
By (A$_1'$) and Hölder's inequality with exponents $r_0=\frac{q_0}{s}$ and $r_0'=\frac{q_0}{q_0-s}$,
\begin{multline*}
    \left( \int_{\mathbb{R}^n} w(x)^s f(x)^s h_2(x) \, dx \right)^{\frac{1}{s}} \le \left( \int_{\mathbb{R}^n} w(x)^s f(x)^s H_2(x)^{\frac{(q - q_0)s}{q_0 (q - s)}} h_2(x)^{\frac{(q_0 - s)q}{q_0 (q - s)}} \, dx \right)^{\frac{1}{s}} \\
    \leq \left( \int_{\mathbb{R}^n} (w(x)^s f(x)^s)^{r_0} H_2(x)^{\frac{(q - q_0)s}{q_0(q - s)}r_0} \, dx \right)^{\frac{1}{r_0 s}} \left( \int_{\mathbb{R}^n} h_2(x)^{\frac{q}{q - s}} \, dx \right)^{\frac{1}{r_0' s}}.
\end{multline*}
By $\|h_2\|_{L^{r'}(\mathbb{R}^n)}=1$, Hölder's inequality with exponents $\frac{q}{q_0}$ and $\big(\frac{q}{q_0}\big)'=\frac{q}{q-q_0}$, and (B$_1'$),
\begin{multline*}
    \|f\|_{L^q(\R^n,w)} 
     \le \left( \int_{\mathbb{R}^n} w(x)^q f(x)^q \, dx \right)^{\frac{1}{q}} \left( \int_{\mathbb{R}^n} H_2(x)^{r'} \, dx \right)^{\frac{q - q_0}{qq_0}} \\
     \leq \| f \|_{L^q(\R^n,w)} \| H_2 \|_{L^{r'}(\mathbb{R}^n)}^{r'\frac{q-q_0}{qq_0}} 
     \leq  2^{r' \frac{q - q_0}{qq_0}} \| f \|_{L^q(\R^n,w)}  < \infty.
\end{multline*}
\medskip

Define $w_0(x) = w(x) H_2^{( \frac{q - q_0}{q - s} ) \frac{1}{q_0}}$; we  claim that  that $w_0 \in A_{p_0, q_0, \B}$. Let $\alpha = \frac{q-s}{q_0-s}$. Then $\alpha >1$ since $q_0<q$. So if we apply Hölder’s inequality with exponents $\alpha$ and $\alpha'=\frac{q-s}{q-q_0}$, the first integral in the definition of $ A_{p_0, q_0, \B}$ becomes
\begin{align*}
    \left( \avgint_{B} \left( w(x) H_2(x)^{\left( \frac{q - q_0}{q - s} \right) \frac{1}{q_0}} \right)^{q_0} dx \right)^{\frac{1}{q_0}} & = \left( \avgint_{B} w^{q_0} H_2(x)^{\frac{q - q_0}{q - s}} dx \right)^{\frac{1}{q_0}} \\
    & =\left( \avgint_{B} w(x)^{q_0} w(x)^{-\frac{s(q-q_0)}{q-s}} (H_2w^s)^{\frac{q - q_0}{q - s}}(x) dx \right)^{\frac{1}{q_0}} \\
    & \leq \left( \avgint_{B} w(x)^q dx \right)^{\frac{1}{q_0 \alpha}} \left( \avgint_{B} (H_2w^s)(x) dx \right)^{\frac{q - q_0}{q_0(q-s)}} \\
    & = \left( \avgint_{B} w(x)^q dx \right)^{\frac{q_0-s}{q_0(q-s)}} \left( \avgint_{B} (H_2w^s)(x) dx \right)^{\frac{q - q_0}{q_0(q-s)}} \\
    & = \left( \avgint_{B} w(x)^q dx \right)^{\frac{p'}{q p'_0}} \left( \avgint_{B} (H_2w^s)(x) dx \right)^{\frac{q - q_0}{q_0(q-s)}} .
\end{align*}
To estimate the second integral in the definition, since $p'=\frac{qs}{q-s}$, $p'_0=\frac{q_0s}{q_0-s}$, by (C$_1'$) we have that
\begin{align*}
    \left( \avgint_{B} w_0(x)^{-p_0'} \, dx \right)^{\frac{1}{p_0'}} & = \left( \avgint_{B} \left(w(x) H_2(x)^{\frac{q - q_0}{q_0 (q - s)}} \right)^{-p_0'}\, dx \right)^{\frac{1}{p_0'}} \\
    & = \left( \avgint_{B} \left(w(x) w(x)^{-\frac{s(q-q_0)}{q_0(q-s)}}  (H_2w^s)^{\frac{q - q_0}{q_0 (q - s)}}(x) \right)^{-p_0'}\, dx \right)^{\frac{1}{p_0'}} \\
    & = \left( \avgint_{B} w(x)^{-\frac{q(q_0-s)}{q_0(q-s)}p'_0}  (H_2w^s)^{-\frac{q - q_0}{q_0 (q - s)}p'_0}(x) \, dx \right)^{\frac{1}{p_0'}} \\
    & \le [H_2w^s]_{\A_{1,\B}}^{\frac{q-q_0}{q_0(q-s)}}\left(\avgint_B (H_2w^s)(x)dx \right)^{-\frac{(q-q_0)}{q_0(q-s)}}\left( \avgint_B w(x)^{-p'} dx \right)^{\frac{1}{p'_0}}.
\end{align*}
If we combine these two inequalities, again by (C$_1'$) we get 
\begin{multline*}
    \left( \avgint_{B} w_0(x)^{q_0} \, dx \right)^{\frac{1}{q_0}} \left( \avgint_{B} w_0(x)^{-p_0'} \, dx \right)^{\frac{1}{p_0'}} \\
    \le [H_2w^s]_{\A_{1,\B}}^{\frac{q-q_0}{q_0(q-s)}} [w]_{\A_{p,q,\B}}^{\frac{p'}{p'_0}} 
    \le  \big( 2\|M_\B \|_{L^{r'}(w^{-s})}\big)^{\frac{q-q_0}{q_0(q-s)}} [w]_{\A_{p,q,\B}}^{\frac{p'}{p'_0}} < \infty.
\end{multline*}

\medskip  
We can now prove inequality \eqref{eqn:extrapol1}.  Since $\|h_2\|_{L^{r'}(\mathbb{R}^n)}=1$, by \eqref{eqn:extrapol0}, Hölder’s inequality with exponents $\frac{p}{p_0}$ and $\big( \frac{p}{p_0}\big)'=\frac{p}{p-p_0}$,  (B$_1'$),and since  $N_{p_0, q_0}$ is increasing, we have that
\begin{align*}
    \| f \|_{L^q(\R^n,w)} & \le \left( \int_{\mathbb{R}^n} w(x)^{q_0} f(x)^{q_0} H_2(x)^{\frac{q - q_0}{q - s}} \, dx \right)^{\frac{1}{q_0}} \| h_2 \|_{L^{r'}(\mathbb{R}^n)}^{\frac{r'}{r'_0s}}\\
    & \le N_{p_0, q_0} ([w_0]_{A_{p_0, q_0, \B}}) \left( \int_{\mathbb{R}^n} (w(x) H_2(x)^{\frac{(q - q_0)}{q_0 (q - s)}} g(x))^{p_0} \, dx \right)^{\frac{1}{p_0}} \\
    & \leq N_{p_0, q_0} ([w_0]_{A_{p_0, q_0, \B}}) \left( \int_{\mathbb{R}^n} (w(x) g(x))^p \, dx \right)^{\frac{1}{p}} \left( \int_{\mathbb{R}^n} H_2^{\frac{(q - q_0)}{q_0 (q - s)} p_0 \frac{p}{p-p_0}} \, dx \right)^{\frac{p - p_0}{p_0 p}} \\
    & = N_{p_0, q_0} ([w_0]_{A_{p_0, q_0, \B}}) \| H_2 \|_{L^{r'}(\mathbb{R}^n)}^{r' \frac{p - p_0}{p_0 p}} \|g\|_{L^p(\R^n,w)}\\
    & \leq 2^{\frac{q - q_0}{q_0 (q - s)}} N_{p_0, q_0} \left( \left( 2\|M_\B \|_{L^{r'}(w^{-s})}\right)^{\frac{q-q_0}{q_0(q-s)}} [w]_{\A_{p,q,\B}}^{\frac{p'}{p'_0}}\right) \|g\|_{L^p(\R^n,w)};
\end{align*}
here we also used that $\frac{q-q_0}{q_0q}=\frac{p-p_0}{p_0p}$ and  $r'=\frac{q}{q-s}$.

\medskip

\textbf{Case IV: $\boldsymbol{p_0=1}$ and $\boldsymbol{q_0<\infty}$}. Fix $1<q_0<\infty$ and let $p,\,q,\,r$ be as in the hypotheses.  In this case,  $s=q_0$, $r_0=1$, and $r=\frac{q}{q_0}>1$. Fix $w\in A_{p,q,\B}$.  We argue as we did in Case III:  first we define the iteration operator $\Rdf'$ and then use duality to get $h_2\in L^{r'}(\mathbb{R}^n), \|h_2\|_{L^{r'}(\mathbb{R}^n)}=1$ such that~\eqref{eqn:dual-caseIII} holds. Define $H_2=\mathcal{R}'(h_2w^s)w^{-s}$. Let $w_0=(H_2w^s)^{\frac{1}{q_0}}$; then property (C$_1'$) implies $w_0^{q_0}\in \A_{1,\B}$, so $w_0\in \A_{1,q_0,\B}$. Moreover, by the definitions, $[w_0^{q_0}]_{\A_{1,\B}}=[w_0]^{q_0}_{\A_{1,q_0,\B}}$. 

By Hölder's inequality and (B$_1'$) 
\begin{multline*}
    \left( \int_{\mathbb{R}^n} w(x)^sf(x)^s h_2(x) dx\right)^{\frac{1}{s}} \le \left( \int_{\mathbb{R}^n} w(x)^s f(x)^s H_2(x)dx \right)^{\frac{1}{s}} \\
    \le \|f\|_{L^q(\R^n,w)} \|H_2\|^{\frac{1}{q_0}}_{L^{r'}(\mathbb{R}^n)} \le 2^{\frac{1}{q_0}} \|f\|_{L^q(\R^n,w)}<\infty .
\end{multline*}
Therefore, to prove~\eqref{eqn:extrapol1}, by \eqref{eqn:extrapol0}, H\"older's inequality, (B$_2'$), and since  $N_{1,q_0}$ is increasing,
\begin{align*}
    \|f\|_{L^q(\R^n,w)} & \le \left(\int_{\mathbb{R}^n} w_0(x)^{q_0} f(x)^{q_0} dx \right)^\frac{1}{q_0} \\
    & \le N_{1,q_0}([w_0]_{\A_{1,q_0,\B}}) \int_{\mathbb{R}^n} w_0(x)g(x)dx \\
    & = N_{1,q_0}([w_0]_{\A_{1,q_0,\B}}) \int_{\mathbb{R}^n} (H_2w^s)^{\frac{1}{q_0}}(x)g(x)dx \\
    & \le N_{1,q_0}([w_0]_{\A_{1,q_0,\B}}) \left( \int_{\mathbb{R}^n} (w(x)g(x))^p \right)^{\frac{1}{p}} \left( \int_{\mathbb{R}^n} H_2(x)^{\frac{p'}{q_0}} dx\right)^{\frac{1}{p'}} \\
    & = N_{1,q_0}([w_0]_{\A_{1,q_0,\B}}) \|H_2\|_{L^{r'}(\mathbb{R}^n)}^{\frac{1}{s}} \|g\|_{L^p(\R^n,w)} \\
    &\le N_{1,q_0}((2 \| M_\B \|_{L^{r'}(w^{-s})})^{\frac{1}{q_0}}) 2^{\frac{1}{q_0}} \|g\|_{L^p(\R^n,w)}.
\end{align*}    
\end{proof}

\medskip

\begin{remark} \label{remark:scalar-matrix}
    We  note a few of the significant differences in the proofs of extrapolation between the scalar and matrix case.  While the boundedness of the scalar iteration algorithm is a straightforward consequence of weighted norm inequalities for the Hardy-Littlewood maximal operator, defining the convex-set valued maximal operator, proving it is bounded on matrix weighted spaces, and proving the existence of the iteration algorithm (Theorem~\ref{thm_7.6}), required developing the theory of convex-set valued functions in~\cite{bownik-cruz-uribe}.   
    
    Moreover, in the matrix case we needed to introduce the auxiliary operator $N_W$ so that the the iteration operators yielded ellipsoid or ball-valued functions.  This allowed us to use Lemma~\ref{lemma:convex-A1} to pass from the convex-set $\A_1^\K$ weights produced by the iteration algorithm  to matrix $\A_{1,\B}$ weights.  
    
    In the scalar case, it is straightforward to use the definition of $\A_{p,q,\B}$ to prove that the weights $w_0$ were in $\A_{p_0,q_0,\B}$.  This follows from the fact that it is easy to prove that if $w_0,\,w_1\in \A_{1,\B}$, then $w_0^{1/p} w_1^{-1/p'}\in \A_{p,\B}$; this is the so-called reverse factorization property of scalar weights.   In the matrix case this is much more difficult to prove; this is hidden in the appeal to Lemma~\ref{matrix-weightsA_r}.

    Finally, we note that in the proof of the scalar case, we can define powers of functions (e.g., $h_1^s$ in Case~I).  But in the matrix case, it requires more care (e.g., $H_1^s$ in Case I of Theorem~\ref{thm:sharp-off-diag-extrapol-pairs}). This is another reason why we introduced the operator $N_W$: if $F$ is a generic convex-set valued function, then $F^s$ is not well defined unless it is ellipsoid-valued.  
\end{remark}

\bibliographystyle{plain}
\bibliography{matrix-extrapol-bases}

\begin{thebibliography}{10}

\bibitem{bownik-cruz-uribe}
M.~Bownik and D.~Cruz-Uribe.
\newblock Extrapolation and factorization of matrix weights.
\newblock {\em Math. Ann.}, 395(3):Paper No. 55, 84, 2026.

\bibitem{buckley1993}
S.~M. Buckley.
\newblock Estimates for operator norms on weighted spaces and reverse jensen
  inequalities.
\newblock {\em Transactions of the American Mathematical Society},
  340(1):253--272, 1993.

\bibitem{preprint-DCU}
D.~Cruz-Uribe.
\newblock Extrapolation and factorization.
\newblock In J.~Lukes and L.~Pick, editors, {\em Function spaces, embeddings
  and extrapolation X, Paseky 2017}, pages 45--92. Matfyzpress, Charles
  University, 2017.
\newblock arXiv:1706.02620.

\bibitem{MR4972891}
D.~Cruz-Uribe and F.~\c{S}irin.
\newblock Off-diagonal matrix extrapolation for {M}uckenhoupt bases.
\newblock {\em Ann. Fenn. Math.}, 50(2):577--609, 2025.

\bibitem{MR2797562}
D.~Cruz-Uribe, J.~M. Martell, and C.~P{\'e}rez.
\newblock {\em Weights, extrapolation and the theory of {R}ubio de {F}rancia},
  volume 215 of {\em Operator Theory: Advances and Applications}.
\newblock Birkh\"auser/Springer Basel AG, Basel, 2011.

\bibitem{MR3544941}
D.~Cruz-Uribe, K.~Moen, and S.~Rodney.
\newblock Matrix {$\mathcal A_p$} weights, degenerate {S}obolev spaces, and
  mappings of finite distortion.
\newblock {\em J. Geom. Anal.}, 26(4):2797--2830, 2016.

\bibitem{10.1016/j.jfa.2010.12.015}
J.~Duoandikoetxea.
\newblock {Extrapolation of weights revisited: New proofs and sharp bounds}.
\newblock {\em J. Func. Anal.}, 260(6):1886 -- 1901, 03 2011.

\bibitem{harboure-macias-segovia88}
E.~Harboure, R.~A. Mac{\'{\i}}as, and C.~Segovia.
\newblock Extrapolation results for classes of weights.
\newblock {\em Amer. J. Math.}, 110(3):383--397, 1988.

\bibitem{hytonenP2010}
T.~Hyt\"{o}nen.
\newblock The sharp weighted bound for general {C}alder\'{o}n-{Z}ygmund
  operators.
\newblock {\em Ann. of Math. (2)}, 175(3):1473--1506, 2012.

\bibitem{MR3204859}
T.~Hyt\"{o}nen.
\newblock The {$A_2$} theorem: remarks and complements.
\newblock In {\em Harmonic analysis and partial differential equations}, volume
  612 of {\em Contemp. Math.}, pages 91--106. Amer. Math. Soc., Providence, RI,
  2014.

\bibitem{MR2652182}
M.~T. Lacey, K.~Moen, C.~P\'erez, and R.~Torres.
\newblock Sharp weighted bounds for fractional integral operators.
\newblock {\em J. Funct. Anal.}, 259(5):1073--1097, 2010.

\bibitem{MR1428988}
F.~Nazarov and S.~Treil.
\newblock The hunt for a {B}ellman function: applications to estimates for
  singular integral operators and to other classical problems of harmonic
  analysis.
\newblock {\em Algebra i Analiz}, 8(5):32--162, 1996.

\bibitem{MR4966577}
Z.~Nieraeth.
\newblock A lattice approach to matrix weights.
\newblock {\em Math. Ann.}, 393(1):993--1072, 2025.

\end{thebibliography}
\end{document}